\documentclass[11pt]{article}
\usepackage{amssymb}
\usepackage{latexsym}
\usepackage{amsmath}
\newtheorem{lemma}{Lemma}
\newtheorem{proposition}{Propostion}
\newtheorem{theorem}{Theorem}
\newtheorem{corollary}{Corollary}
\newtheorem{claim}{Claim}
\newenvironment{proof}[1][Proof]
           {\medbreak\noindent {\em #1. \enspace}}
           {\enspace q.e.d. \par \medbreak}
\newenvironment{definition}[1][Definition]
           {\medbreak\noindent {\em #1. \enspace}}
           {\par \medbreak}
\newenvironment{remark}[1][Remark]
           {\medbreak\noindent {\em #1. \enspace}}
           {\par \medbreak}
\newenvironment{remarks}[1][Remarks]
           {\medbreak\noindent {\em #1. \enspace}} {\par \medbreak}
\makeatletter
\@addtoreset{equation}{section}
\makeatother

\def\AA{{\mathbb A}}

\def\RR{{\fontsize{12}{14pt}\selectfont\mathbb{R}}}

\newcommand{\fg}{\mathfrak{g}}

\newcommand{\s}{\sigma}
\renewcommand{\o}{\omega}
\newcommand{\al}{\alpha}
\newcommand{\be}{\beta}

\renewcommand{\d}{\delta}
\newcommand{\e}{\varepsilon}
\newcommand{\g}{\gamma}
\newcommand{\la}{\lambda}
\newcommand{\La}{\Lambda}

\renewcommand{\t}{\tau}
\renewcommand{\th}{\theta}

\newcommand{\z}{\zeta}

\newcommand{\sm}{\setminus}

\newcommand{\pt}{\partial}
\newcommand{\ra}{\rightarrow}
\newcommand{\lra}{\longrightarrow}
\newcommand{\dra}{\Rightarrow}

\renewcommand{\rm}{\textrm}

\newcommand{\cA}{\mathcal{A}}

\newcommand{\cE}{\mathcal{E}}
\newcommand{\cM}{\mathcal{M}}

\newcommand{\cS}{\mathcal{S}}

\newcommand{\cL}{\mathcal{L}}

\newcommand{\cC}{\mathcal{C}}
\def\End{\mathop{\rm{End}}\nolimits}

\def\Ker{\mathop{\rm{Ker}}\nolimits}
\def\Im{\mathop{\rm{Im}}\nolimits}
\def\Stab{\mathop{\rm{Stab}}\nolimits}
\def\GL{\mathop{\rm{GL}}\nolimits}
\def\Ad {\mathop{\rm{Ad}}\nolimits}
\def\YM {\mathop{\rm{YM}}\nolimits}
\def\injrad {\mathop{\rm{injrad}}\nolimits}

\def\Id {\mathop{\rm{Id}}\nolimits}
\def\div {\mathop{\rm{div}}\nolimits}
\def\grad {\mathop{\rm{grad}}\nolimits}
\def\Aut {\mathop{\rm{Aut}}\nolimits}

\def\Re {\mathop{\rm{Re}}\nolimits}
\def\dist {\mathop{\rm{dist}}\nolimits}

\author{BAOZHONG YANG \footnote{Revised on Feb. 28, 2002.}}
\title{The uniqueness of tangent cones for Yang-Mills connections with isolated singularities}
\date{}
\begin{document}
\maketitle

\begin{abstract}
 We proved a uniqueness theorem of tangent connections for a
 Yang-Mills connection with an isolated singularity with a quadratic
 growth of the curvature at the singularity. We also obtained control
 over the rate of the asymptotic convergence of the connection to the
 tangent connection. The rate of convergence can be strengthened when
 the tangent cone is integrable. There are parallel results for the
 cones at infinity of a Yang-Mills connection on an asymptotically
 flat manifold. We also gave an application of our methods to the
 Yang-Mills flow and proved that the Yang-Mills flow exists for all
 time and has asymptotic limit if the initial value is close to a
 smooth local minimizer of the Yang-Mills functional.
\end{abstract}

%
%
\section*{Introduction}
The studies of singularities have been important in geometric
analysis. Although in many cases our prime interests are in studying
smooth solutions of certain geometric PDE's, singular solutions arise
naturally as the limits of sequences of smooth ones. The local
behavior of a singularity is largely determined by the so called
tangent cones. Tangent cones have been important in various
kinds of regularity theories, and the use of them has been especially
successful in the theory of minimal surfaces and harmonic maps (see
for example the books \cite{SL2}, \cite{SL4}, \cite{G} and the paper
\cite{ScU}). 

The uniqueness of tangent cones at a singular point of a geometric
object is equivalent to the asymptotic convergence of the object to a
tangent cone. This case provides the simplest possible nontrivial
singularity and the local singularity of the geometric object would be
better understood if we also know the rate of its asymptotic
convergence to the unique cone.

The important work by Leon Simon \cite{SL1} gave a general asymptotic
convergence theorem for the solutions to a class of nonlinear
evolution equations which arise naturally in certain geometric
variational problems, especially in problems associated with tangent
cones. In the same paper, L. Simon applied the general theorem to
show the uniqueness of tangent objects for minimal surfaces and energy
minimizing maps with the assumption that there exists a tangent object
with an isolated singularity. An earlier work of W. Allard and
F. Almgren \cite{AA} proves the case of minimal surfaces with the
integrability assumption on the tangent cone. 

This result and the methods of Leon Simon and its variants have been
used and studied by many other people in different settings. The work
Cheeger and Tian \cite{CT} proved a uniqueness theorem for cone
structures at infinity of Ricci flat manifolds assuming that the cone
at infinity being integrable. The work Morgan, Mrowka and Ruberman
\cite{MMR} studied the asymptotic limits of anti-self-dual connections
with finite energy on a cylindrical 4-manifolds and gave applications
to computations of Donaldson invariants.

In this paper, we proved the uniqueness of tangent cones of a
Yang-Mills connection with an isolated singularity on manifolds of
dimension greater than $4$, under a quadratic growth assumption on the
curvature of the connection (see Theorem~\ref{th:1.1}). We estimated
the rate of asymptotic convergence of the connection to its cone. On
the other hand, we have a faster convergence rate would be fast if we
assume the tangent cone is integrable. There are parallel results for
the cones at infinity of a Yang-Mills connection on an asymptotically
flat manifold. The author has constructed examples of Yang-Mills
connections with given tangent cones and with different rates of
convergence to the cones in \cite{Y}.

The difficulties of our proof mainly lies in the degenerate elliptic
nature of the Yang-Mills equation. Under a good transverse gauge, the
Yang-Mills equation become elliptic. If we have long time existence of
the gauge, then our result will follow from an application of
L. Simon's result. The long time existence of this gauge, however,
depends on the asymptotic convergence which we want to show. Our
strategy to solve this dilemma is to find a suitably defined gauge and
prove the long time existence of the gauge and the solution in one
shot. In the proof of long-time existence and convergence in Section
3, we modified L. Simon's methods in \cite{SL1} to our case. The
method is to divide the existence interval of the solution into three
parts on which the norms of the solution have different growth
behaviors. Roughly speaking, the norm of the solution is exponentially
decreasing on the first interval, it is changing slowly on the second
interval, and it is exponentially increasing on the third
interval. These behaviors are modelled on those of the solutions to
the linearized equation. One can then control the norms of the
solution on the three intervals using different techniques. The first
interval is easy. For the second interval we used the variational
inequality by L. Simon (see Section 3.3) and for the last interval we
have to use the properties of the gauge we constructed in Section 2.3,
especially the property that under this gauge the time derivative of
the connection is uniformly small on the existence interval. 

One might expect that there are similar results for tangent cones for
more general (not necessarily isolated) singularities, say,
singularities which are higher dimensional submanifolds or
subvarieties. The methods of this paper, however, do not directly
apply to the more general case.

An application of the gauge we used gives a result for the Yang-Mills
flow. We showed that the Yang-Mills flow exists for all time and has an
asymptotic limit if the initial value is close to a smooth local
minimizer of the Yang-Mills functional (see Corollary~\ref{co:4.2}).

This paper is organized into four parts. In Section 1 we review
backgrounds, set up notations and state our main results. In Section
2, we describe our procedure of fixing a suitable transverse gauge and
the related estimates for connections under this gauge. In Section 3
we prove an asymptotic convergence theorem of certain solutions of a
class of evolution equations and this gives long time existence of the
gauge constructed in Section 2 and convergence of the connection at
once. We also showed how we control the rate of the asymptotic
convergence. Finally in Section 4, we proved the result stated above
for the Yang-Mills flow.

{\bf Acknowledgment} This work is based on the author's thesis at
Massachusetts Institute of Technology. I am indebted to my advisor
Gang Tian for his support and guidance in this work. I'd like to thank
David Jerison, Richard Melrose, Tomasz Mrowka, Rich Schoen and Leon
Simon for helpful suggestions and discussions. I also want to thank
people at Brown, Columbia and Princeton who have attended my talk
there and given helpful suggestions and encouragement.

%
%
\section{Background, Notations and Main Results}
\label{sec:1}
%
\subsection{Background and main theorems}
\label{sec:1.1}
We shall give a very brief background here. One may find more detailed
information about Yang-Mills connections in, for example, the
excellent book \cite{DK} and some recent analytic results we used in
\cite{T}. We assume that $M$ is an n-dimensional manifold and $P$ is a
principal bundle on $M$ with compact structure group $G$ which has Lie
algebra $\fg$. Suppose $E$ is a vector bundle associated to $P$ with a
 linear representation $\rho: G \ra \GL(V)$, where $V$ is the
fiber type of $E$. $\Aut P = P \times_{\Ad} G$ is the principal bundle
associated to the $\Ad$ representation and we denote by $\fg_E$ the
associated bundle to $\Aut P$ with the differential of $\rho$, $d\rho:
\fg \ra \End(V)$.  $\fg_E $ is a subbundle of $\End(E)$. The
connections considered in this paper will always be $G$-connections.

Let $A$ be a connection on $E$ which corresponds to a covariant
derivative $\nabla_A: \Gamma(E) \ra \Omega^1(E)$, which is
$G$-equivariant . Under a trivialization $\phi: E|_{U} \ra U
\times V$ on an coordinate open subset $U \in M$, $A$ has local
expression $ \nabla_A = d+A_{\phi} $, where $A_{\phi} \in
\Omega^1(\fg_E|_U) \subset \Omega^1(\End(E|_U))$
. We often suppress the subscript $\phi$ and use $A$ for
this local expression when the trivialization is clear in the
context. The covariant differentiation $\nabla_A$ extends naturally
to sections of various tensor bundles on $M$ with values in $E$ or $\End (E)$.
$d_A = \wedge \circ \nabla_A$ gives the coupled exterior differential
on forms on $M$ with values in $E$ or $\End (E)$.
In particular, in local coordinates,
\begin{align}
 & d_A (\al) = d \al + A \wedge \al, \\
 & d_A (\be) = d \be + [A, \be],
\end{align}
for $\al \in \Omega^*(E)$ and $\be \in \Omega^*(\End(E)) $.

For a connection $A$ on $E$, the operator $d_A \circ d_A: \Omega^*(E)
\ra \Omega^{*+2}(E)$ is given by the algebraic operator $F_A
\wedge: \Omega^*(E) \ra \Omega^{*+2}(E)$, where $F_A \in \Omega^2(\End(E))$ is the curvature of $A$ and
locally has the expression
\begin{equation}
F_A = dA + A\wedge A = dA + \frac{1}{2} [A, \, A].
\end{equation}
The curvature $F_A$ satisfies the Bianchi identity
\begin{equation}
\label{eq:1.1.1}
d_A F_A = 0.
\end{equation}

Let $\cA = \{ \text{all connections  on }E \}$. The gauge group,
$\Gamma(\Aut P)$, consists of bundle automorphisms of $E$ that
preserve the base. The gauge action of $\Gamma(\Aut P)$ on $\cA$ is given by
$$ \nabla_{g(A)}(v) = g \circ \nabla_A \circ g^{-1} (v), $$
for $g \in \Gamma(\Aut P)$, $A \in \cA$ and $v \in \Gamma(E)$.
(Another convention in literature is to let $\nabla_{g(A)} = g^{-1}
\circ \nabla_A \circ g$.)  Under local trivialization, if we write
$\nabla_A = d + A$, $\nabla_{g(A)} = d + g(A)$, then
\begin{equation}
g(A) = g A g^{-1} - dg\, g^{-1}.
\end{equation}
We have, as the curvature is a tensor,
\begin{equation}
F_{g(A)} = g F_A g^{-1}.
\end{equation}

Now we assume that $M$ is oriented and is given a Riemmanian structure
$g$ and assume that $E$ is given a $G$-invariant metric $h$. The
Yang-Mills functional $\YM: \cA\rightarrow \RR$ is defined by
$$\YM(A)= \int_M\left| F_A\right|^2 dV_g.$$
Here $|\, \cdot\,|^2$ is given by the metrics $g$ and $h$.

We call $A$ a {\em Yang-Mills connection} if and only if $A$ is a
critical point of $\YM$ on $\cA$. In other words, $A$ is Yang-Mills if
and only if for any continuously
differentiable family of connections $\{A_t\}_{-\e <t <\e},$
$$ \left.  \frac{d}{dt}\right|_{t=0} \YM(A_t) = 0. $$
A Yang-Mills connection $A$ satisfies the following Yang-Mills
equation, which is the Euler-Lagrange equation for $\YM$,
\begin{equation}
\label{eq:1.1}
d_A^{*}F_A=0.
\end{equation}
In (\ref{eq:1.1}), $d_A^{*}: \Omega^{*}(\End(E)) \ra
\Omega^{*-1}(\End(E))$ is the formal adjoint of $d_A$. 
The Bianchi identity (\ref{eq:1.1.1}) and Yang-Mills equation
(\ref{eq:1.1}) together give a system of equations which is a
nonlinear analogue of the equation for harmonic forms.  

It follows by direct calculation that under local coordinates, if $\al \in \Omega^p(\End(E))$, then
$$d_A^{*}(\al) = d^* \al +(-1)^{(p-1)n +1} *[A, *\al] = (-1)^{(p-1)n +1}*d_A*
\al.$$

By using variations generated by a vector field on $M$, we have the
following {\em first variation formula}  for {\em smooth} Yang-Mills connections
\begin{equation}
\label{eq:var}
\int_M |F_A|^2 \div X - 4 \sum_{1\leq i < j \leq n}\langle F_A
( \nabla_{e_i}X, e_j), \,F_A(e_i, e_j) \rangle dV_g =0.
\end{equation}
 This formula is true for any compactly supported
$C^1$ vector field $X$ on $M$.  \\[1.5ex]
\noindent {\em Definition. } A connection $A$ (possibly with
singularities) is {\em stationary} if the first variation formula
\eqref{eq:var} is true for $A$ with any compactly supported $C^1$
vector field $X$ on $M$.\\

Assume that $A$ is a stationary connection. By using a cutoff of the
radial vector field $X = \sum_{i=1}^{n} x_i \frac{\pt}{\pt x_i} $, we
have the following important monotonicity formula by Price \cite{P}
(see \cite{T} for a variant version). Let $\injrad(x)$ denote the
injective radius of $x \in M$.
\begin{proposition}
\label{prop:1.2.1} 
If $A$ is a stationary connection on $E$, then for any $x \in M$,
there exists positive constants $\Lambda = \Lambda(x)$ and $r_x <
\injrad(x)$ which only depend the supremum bound of curvature of $M$,
such that if $0 < \sigma < \rho \leq r_x$, then
\begin{align}
\label{eq:m}
\rho ^{4-n}e^{\Lambda \rho}& \int_{B_\rho (x)}\left| F_A\right| ^2dV_g-\sigma
^{4-n}e^{\Lambda \sigma}\int_{B_\sigma (x)}\left| F_A\right| ^2dV_g\\
& \geq 4\int_{B_\rho (x)\setminus B_\sigma \left( x\right)
}r^{4-n}\left| \frac \partial {\partial r} \rfloor F_A \right|^2 dV_g
\nonumber
\end{align}
\end{proposition}

\begin{remark}
We remark here that if $M \cong \RR^n$ is the Euclidean space, then we
can take $\Lambda = 0$ and $r_x = \infty$ in the above proposition, and
furthermore, the equality in \eqref{eq:m} holds. In
this case, we have, by differentiating the equation \eqref{eq:m}
\begin{equation}
\label{eq:m1}
\rho^{5-n} \int_{\pt B_{\rho}(x)} |F_A|^2 -  (n-4) \rho^{4-n}
\int_{B_{\rho}(x)} |F_A|^2 = 4 \rho^{4-n} \int_{\pt
B_{\rho}(x)} \left| \frac \partial {\partial r} \rfloor F_A \right|^2
\geq 0.
\end{equation}
\end{remark}



Using the compactness theorems for $L^p$ connections with $p > n$ in
Uhlenbeck \cite{U2}, it is easy to show the following proposition.
\begin{proposition}
\label{prop:1.3.2}
Assume $\{ A_i\}$ is a sequence of smooth Yang-Mills connections on a
bundle $(E, h_i)$ over $(M, g_i)$, with $ |F_{A_i}|$ locally uniformly
bounded in $i$. Assume that $g_i$ and $h_i$ converges smoothly to
metrics $g$ and $h$ on compact sets of $M$. Then there exist a
subsequence $\{A_{i_j}\}$, gauge transformations $\sigma_j $ on $M$
and a smooth Yang-Mills connection $A$ on $E$ over $M$ such that on
any compact set $K\subset M$, $\sigma_j(A_{i_j})$ converges to $A$ in
$C^\infty $ topology. 
\end{proposition}

For more general compactness theorems with only uniform bounds on the
$L^2$ norms of the curvatures, we need to use a priori estimates to
bound the pointwise norm of curvature and the convergence will be
possible only away from a blow-up set -- where the curvature energy of
the sequence concentrate. This blow-up and weak compactness is
familiar in other settings, as in the theories of minimal surfaces and
harmonic maps. We refer the reader to Nakajima \cite{Na} and Tian
\cite{T} for some general compactness theorems.

We now assume that $\dim M \geq 5$, $x_0 \in M$ and $A$ is a
Yang-Mills connection on a bundle $E$ over $M \setminus \{ x_0 \}$,
i.e., $A$ has an isolated singularity at $x_0$. By a standard cutoff
argument, we may show that \eqref{eq:var} is true for $A$ with any
$C^1$ compactly supported vector field $X$ and $A$ is stationary on
$M$. We also assume that the curvature of $A$ satisfies the following
quadratic growth condition in a neighborhood $U$ of $x_0$:
\begin{equation}
|F_A (x) | \leq C r^{-2}, \quad \forall x \in U,
\end{equation}
where $r= \dist(x, x_0)$ on $M$. We may identify $U$ and an open set,
say, $B_1(0)$ in $T_{x_0}M \cong \RR^n$, and denote the induced
pullback metric on $B_1(0)$ also by $g$.  For any $\la \in (0,1) $,
let $\tau_\lambda : T_x M \ra T_x M $ be the scaling map given by
$\tau_\la : v \mapsto \la v$. We may also lift $\tau_{\lambda}$ to
a map on the bundle $E$ via a local trivialization. Define the scaling
of the metrics $g$, $h$ and the connection$A$ by
\begin{equation}
\label{eq:1.4.1}
g_{\la} = \la^{-2} \t_{\la}^*g, \quad h_{\la} =  \t_{\la}^*h, \quad
A_{\la} = \t_{\la}^* A
\end{equation}
Then for any $0< \la <1$, $A_{\la}$ is a stationary Yang-Mills
connection with respect to $g_{\la}$ with an isolated singularity at
$0$. Since $A$ is defined on $B_1(0) \setminus \{ 0 \}$, it follows
that $A_{\al}$ is defined on $B_{\la^{-1}}(0) \setminus \{ 0\}$. The
metrics $g_i$ and $h_i$ converges smoothly to the standard metrics $g$
on $\RR^n$ and $h$ on a trivial bundle coming from the metrics at $0$.
The compactness property Prop.~\ref{prop:1.3.2} now implies that for
any sequence $\la_i \ra 0$, there is a subsequence that after gauge
transformations, converges on compact sets of $\RR^n \setminus \{ 0
\}$ to a smooth Yang-Mills connection $A_0$ on $\RR^n \setminus \{ 0
\}$.

\begin{definition}
$A_0$ as above is called a {\em tangent connection} of the 
connection $A$ at $x_0$.
\end{definition}
  
As in the case of minimal surfaces and harmonic maps, we expect $A_0$
to be radially homogeneous, i.e., a cone. Indeed, we have the following
property for a tangent connection at a point (\cite[5.3.1]{T}).
\begin{lemma}
\label{lem:1.4.1} 
With notations as above and let $r = \dist(x, 0)$, we have
\begin{equation}
\label{eq:1.4.3} 
\frac{\pt}{\pt r} \rfloor F_{A_0}(x) = 0, \quad \forall x \in \RR^n
\setminus \{ 0 \}.
\end{equation}
\end{lemma}
\begin{proof}
The stationarity of $A_0$ follows from the stationarity of
$\t_{\la_i}^*A$ and the strong convergence on compact sets of $\RR^n
\setminus \{ 0 \}$, or follows directly since singular set of $A_0$ is
isolated, with codimension at least 5. Since $\RR^n$ has the flat
metric, by the remark following Prop.~\ref{prop:1.2.1}, we may take $\La = 0$ in the monotonicity formula (\ref{eq:m}) and
\begin{align}
\label{eq:1.4.5}
\rho ^{4-n}& \int_{B_\rho (0)}\left| F_{A_0}\right| ^2dx-\sigma
^{4-n}\int_{B_\sigma (0)}\left| F_{A_0}\right| ^2 dx\\
& = 4\int_{B_\rho (0)\setminus B_\sigma (0)} r^{4-n}\left| \frac \partial {\partial r} \rfloor F_{A_0} \right|^2 dx
\nonumber
\end{align}
However, for any $\rho >0$, we have
\begin{align*}
\rho^{4-n} \int_{B_\rho(0)}& |F_{A_0}|^2 dx 
 = \lim_{i \ra \infty} (\la_i \rho)^{4-n} \int_{B_{\la_i\rho}(x)}
|F_{A}|^2 dV_g \\
& = \lim_{r \ra 0} r^{4-n} \int_{B_{r}(x)} |F_{A}|^2 dV_g \geq 0
\end{align*}
Therefore, both sides of (\ref{eq:1.4.5}) are zero and
\begin{equation}
\label{eq:1.4.6} 
\int_{B_\rho (0)\setminus B_\sigma (0)} r^{4-n}\left| \frac \partial
{\partial r} \rfloor F_{A_0} \right|^2 dx = 0
\end{equation}
This implies (\ref{eq:1.4.3}).
\end{proof}

With the conclusion of the lemma, it is a standard fact that after a
smooth gauge transformation on $\RR^n \setminus \{ 0 \}$, $A_0$ can be
made into radially homogeneous, i.e., $A_0 = p^*(A_0')$ for some
smooth Yang-Mills connection $A_0'$ on the trivial bundle over
$S^{n-1}$, where $p: \RR^n \setminus \{0\} \ra S^{n-1}$ is the natural
projection. We shall assume later that all tangent connections are in
this radially homogeneous form.  With some abuse of terminology, we
shall call this $A_0'$ a {\em tangent connection }or a {\em
tangent cone} of $A$ at $x_0$.

\begin{remark}
Our definition of tangent connections are of course not the most
general. In fact, using the general compactness theorems in \cite{Na}
and \cite{T} which only assumes a uniform bound on the $L^2$ norms of
the curvature and the monotonicity formula, a tangent connection may
be defined for a stationary Yang-Mills connection with a singularity
set $S$ at a singular point $x_0 \in S$. We need some additional
assumptions on $S$, though, for example, we may have to assume that
$C(S) = \{ x = x_0 + t (y - x_0): y \in S, t \in (0, \infty) \}$ has
zero $(n-4)$-dimensional Hausdorff measure, or assume that $S$ is
stratified by manifolds with dimension not greater than $n-4$. These
assumptions on $S$ are necessary because we don't have the best
possible removable singularity theorems for stationary Yang-Mills
connections (for a good one, see \cite{TT}). In the case of
Hermitian-Yang-Mills connections, we only need to assume the
$(n-4)$-dimensional Hausdorff measure of $S$ is finite (see
\cite{TY}).
\end{remark}

Now we are ready to state our main theorems. We assume that $(M, g)$
is of dimension $n \geq 5$, $ x_0 \in M$ and $(E, h)$ is a bundle over
$M \setminus \{ x_0 \}$.
\begin{theorem}
\label{th:1.1}
Let $A$ be a smooth Yang-Mills connection on the bundle $E$
over $M\sm \{ x_0 \}$. Assume that there exists $C >0$ and a neighborhood
$U$ of $x$ such that 
$$ | F_A(x)| \leq C r^{-2}, \forall x \in U, $$
where $r = \dist(x,
x_0)$. Then the tangent connection of $\AA$ is unique up to gauge
transformations. In other words, there exists a smooth Yang-Mills
connection on $E|_{S^{n-1}(x_0)}$ over the unit sphere around $x_0$ and
a gauge transformation $\tau$ on $E|_{M \setminus \{x_0\}}$ such that
$$ A(r) \stackrel{C^{\infty}}{\lra} A_0 \quad
\rm{as } r \ra 0,$$
where $A(r) = \tau_r^*(A)|_{S^{n-1}(x_0)}$ is the rescaled connection
under the scaling map $\tau_r: B_{r^{-1}}(x_0) \ra B_1(x_0)$ as
before. Furthermore, there exists constants $C_k > 0$ and $\al >0$
depending on $A$, such that
$$ | A(r) - A_0 |_{C^k(S^{n-1})} \leq C_k | \log r|^{-\al}. $$
\end{theorem}
We call a smooth Yang-Mills connection $A_0$ on $S^{n-1}$ {\em integrable}
if for every solution $a \in \Omega^1(\fg_E)$ of
\begin{equation}
\label{eq:3.5.17}
L a = \Delta_{A_0} a + (-1)^n *[a, *F_{A_0}] = 0,
\end{equation}
where $L$ is, as in Section 3.1, the linearization of $d_{A_0
+a}^*F_{A_0 + a}$ at $0$, there exists a path of Yang-Mills
connections $A(t)$, $t \in (-\e, \e)$ with $A(0) = A_0$ such that
\begin{equation}
\label{eq:3.5.18}
\left.\frac{\partial }{\partial t}\right|_{t=0} A(t) = a.
\end{equation}
This has the geometric meaning that $A_0$ has an integrable
neighborhood in the moduli space of smooth Yang-Mills connections on
$S^{n-1}$ with tangent space at $A_0$ being given by the Jacobi fields
at $A_0$, i.e., the solutions to (\ref{eq:3.5.17}).

\begin{theorem}
\label{th:1.2}
If in addition to the hypotheses of theorem 1, we assume that the
tangent connection $A_0$ of $A$ is integrable, then there exists $\al >
0$ and $C_k > 0$ such that we have the following better control of convergence,
$$ |A(r) - A_0|_{C^k(S^{n-1})} \leq C_k r^{-\al}. $$
\end{theorem}

%
%
%
\begin{remark}
1) We remark that in the proof of Theorem~\ref{th:1.1}, the
monotonicity formula is only used to prove the desired rate of
convergence. In \cite{Y}, the author was able to construct examples of
connections with a given smooth tangent cone with the convergence
rates as in Theorem~\ref{th:1.1} and Theorem~\ref{th:1.2}
respectively. These examples are defined on a ball with certain
boundary values and the construction is achieved by a
deformation-perturbation method. It would be interesting to construct
examples defined on compact manifolds.\\ 2) There are parallel results
of Theorem~\ref{th:1.1} and \ref{th:1.2} for the tangent cones at
infinity of a Yang-Mills connection on bundles over $\RR^n$ ($n \geq
5$), or more generally on an asymptotically flat manifold. We need to
assume quadratic decay of the curvature there. The proofs will be the
same. We leave the formulation to the reader.
\end{remark}
%




\subsection{Cylindrical coordinates and notations}
\label{sec:2.1} 
Since our problem in Theorem~\ref{th:1.1} is of local nature, we may
assume that $E$ is a bundle over the disk $B_2(0) \setminus \{ 0 \}$. Consider
the cylindrical coordinates, $\phi : B_{r_0}(0)\ra S^{n-1} \times
[t_0, \infty), t_0 = - \log r_0$ defined by
$$\phi (x) = ( \o (x), t(x)) = (\frac{x}{|x|}, -\log|x|).$$ We may
assume that the disk $B_2(0)$ has the standard metric and the bundle
$E$ has the product metric, because nonstandard metrics will only give a
perturbation which is exponentially decaying in $t$ as $t \ra \infty$
in the cylindrical coordinates and will not affect our proof later. We
identify $E$ with $(\phi^*)^{-1}(E)$.

Assume that $\tilde{A}$ is a Yang-Mills connection on the bundle over
$B_2(0) \setminus \{ 0 \}$ and assume
$$ (\phi^*)^{-1} (\tilde{A}) = A(t) + \be(t)dt,$$
where $A \in \Omega^1(\End E)$, $\be \in \Omega^0(\End E)$.

\begin{lemma}
\label{lem:2.1.1} 
The Yang-Mills equation $d_{\tilde{A}}^{*} F_{\tilde{A}} =0$ is
equivalent to the following system of equations,
\begin{eqnarray}
\label{eq:2.1}
&\ddot{A} - (n-4)\dot{A} - d_A^* F_A -  d_A \dot{\be} + (n-4) d_A \be + (-1)^{n+1}\ast[\be,
\ast d_A \be]=0&\\
\label{eq:2.2}
&d_A^* (\dot{A} - d_A \be) =0&
\end{eqnarray}
\end{lemma}
\begin{proof}
Recall that $ d_A^*(\xi) = (-1)^{n(p-1)+1}( \ast d \ast \xi + \ast
[A,\ast \xi] )$, for $\xi \in \Omega^p(\End E)$ over an
$n$-dimensional manifold. In considering $(\phi^{*})^{-1} \circ
d^*_{\tilde{A}}$, we know that $d$ commutes with $(\phi^*)^{-1}$, therefore
we need to consider the behavior of $*$ under $(\phi^*)^{-1}$. We
denote by $g_0$ the standard metric in $\RR^n$, by $g$ the standard
product metric on $S^{n-1} \times \RR$, and by $\tilde{g} =
(\phi^*)^{-1} (g_0)$ the pushforward metric of $g_0$. Denote by
$*_{n0}$, $*_n$, $\tilde{*}_n$ and $*$ the Hodge operator associated
to respectively $g_0$, $g$, $\tilde{g}$ and the standard metric
$S^{n-1}$. Let $d_{n0}$, $d_n$ and $d$ be the exterior differential on
respectively $\RR^n$, $S^{n-1} \times \RR$ and $S^{n-1}$. If $\{d\o^1,
\ldots, d\o^{n-1}, dt\}$ is a local orthonormal basis of $T^1
(S^{n-1}\times \RR)$ with respect to $g$, by definition of $\phi$, $\{
e^{-t}d\o^i, e^{-t}dt \}$ constitutes an orthonormal basis for
$\tilde{g}$. Hence
\begin{equation}
\label{eq:2.1.1} 
\tilde{g} = e^{2t} g, ~~\langle \al, \be \rangle_{\tilde{g}} = e^{2 \deg
\be t} \langle \al, \be
\rangle_{g},~~ dV_{\tilde{g}} = e^{-nt}dV_{g}.
\end{equation}
Therefore
$$ \al \wedge \tilde{*}_n \be = \langle \al, \be \rangle_{\tilde{g}}
dV_{\tilde{g}} = e^{2 \deg \be t} \langle \al, \be \rangle_{g}
e^{-nt}dV_{g} = e^{-(n- 2 \deg \be)t}  \al \wedge *_n \be. $$
That is
$$ \tilde{*}_n \be = e^{-(n- 2 \deg \be)t} *_n \be. $$
We can now carry out the calculation, first
$$ d_{\tilde{A}}^{*} F_{\tilde{A}} = (-1)^{n+1} (
*_{n0}d_{n0}*_{n0}F_{\tilde{A}} + *_{n0}[\tilde{A},
*_{n0}F_{\tilde{A}}]) $$ and we have
$$ (\phi^*)^{-1}(F_{\tilde{A}}) = F_A - ( \dot{A} - d_A \be) dt $$
For
simplicity of notation, let $\eta = \dot{A} - d_A $.
\begin{align*}
(\phi^*)^{-1}(*_{n0}& d_{n0}*_{n0}F_{\tilde{A}})= \tilde{*}_n d_n
\tilde{*}_n( F_A - \eta dt ) \\
& = e^{(n-2)t}*_n d_n e^{-(n-4)t} *_n
(F_A - \eta dt) \\
& = (-1)^ne^{2t}\{ (\dot{\eta} - (n-4) \eta - d^*F_A + 
(*d*\eta) \wedge dt \} \\
(\phi^*)^{-1}(*_{n0} & [\tilde{A},
*_{n0}F_{\tilde{A}}]) = \tilde{*}_n [A+\be dt, \tilde{*}_n(F_A - \eta
dt)] \\
&= e^{2t}*_n[A+\be dt, *_n (F_A - \eta dt)]\\
&= e^{2t}\{  - *[A, *F_A] + *[\be, *\eta])+ (-1)^{n-1} *[A, *\eta]dt\}
\end{align*}
Combining the above two equalities, we obtain
$$ (\phi^*)^{-1}(d_{\tilde{A}}^{*} F_{\tilde{A}}) = - e^{2t} \{
 - (\dot{\eta} - (n-4)\eta - d^*_A F_A
+(-1)^n *[\be, *\eta])+ d^*_A \eta \wedge dt \} $$
$\tilde{A}$ is Yang-Mills if and only if the $dt$ part and the part
 without $dt$ above are zero. This together with the observation that
 $[\dot{A}, \be] = (-1)^n *[\be, * \dot{A}]$ give the desired system
 (\ref{eq:2.1}) and (\ref{eq:2.2}).
\end{proof}

We shall use the above cylindrical coordinates in our proof and
consider solutions to the system (\ref{eq:2.1}) and
(\ref{eq:2.2}). Our arguments in this and the next chapter will work
for more general settings, so we assume in the following of this paper
that $(E, h)$ is a Euclidean vector bundle on a compact $(n-1)$
-dimensional Riemannian manifold $(M, g)$ $(n \geq 5)$ with a
compact structure group $G$, and $P$ is the associated principal
bundle.  We let the bundle $E \times [t_0, \infty)$ and the manifold
$M \times [t_0, \infty)$ have the product metrics.

We shall use letters with tildes, $\tilde{A}$, $\tilde{A}_1$, etc to
represent connections on the bundle $E\times I$ on $M\times I$, where
$I$ is an interval of possibly infinite length. And we use $A$, $A_0$,
$B$, etc to represent connections on bundle $E$ on $M$.

For simplicity, we shall assume that $E$ is trivial. For nontrivial
$E$, our arguments still work if we fix a smooth connection $A_0$ on
$E$ and replace $A$ by $A - A_0$, $d$ and $d^*$ by $d_{A_0}$ and
$d_{A_0}^*$ etc.

\begin{definition}
A connection $\tilde{A}$ on the bundle $E
\times [t_0, \infty)$ a {\em Yang-Mills} connection if it satisfies
the system (\ref{eq:2.1}) and (\ref{eq:2.2}).
\end{definition}

Fix $\mu \in (0, 1)$ and let $k$ be a nonnegative integer.  We shall
use in the following various H\"older norms and spaces. The norms may be
defined by using partition of unity arguments, for example. The reason
we use H\"older norms rather than Sobolev norms is that the
restriction to submanifolds of $C^{k, \mu}$ functions are $C^{k, \mu}$
functions, which makes many statements simpler. However, it should be
observed that the Sobolev norms also work for the proof with
appropriate attention to the corresponding trace properties of Sobolev
functions.

Let $A(t)$ be a section of a bundle over $M \times I$, with $I$ an
interval. We shall use the following abbreviation,
\begin{eqnarray*}
|A(t)|_{C^{k, \mu}} &=& |A(t)|_{C^{k, \mu}(M)},\\ 
|A(t)|_{C^{k, \mu}(I)} &=& |A(t)|_{C^{k, \mu}(M \times I)}.
\end{eqnarray*}
Note that the $| \cdot |_{C^{k, \mu}}$ norm doesn't count in the time
derivatives. We define
\begin{multline}
\label{eq:2.1.7}
| A(t) |_{C^{(k,l), \mu}(I)}  =
        \sum_{
                0 \leq j \leq k,
                0 \leq i+j \leq k+l
              }
\sup_I |\nabla_M^i \frac{\partial^j}{\partial t^j} A(t) |_{C^{0}} \\+ \sum_{
                0 \leq j \leq k,
                0 \leq i+j = k+l
              }
\sup_I |\nabla_M^i \frac{\partial^j}{\partial t^j} A(t) |_{C^{0, \mu}}
\end{multline}
The spaces $C^{(k, l), \mu}(I)$ are defined by those bundle sections
such that the norms as in (\ref{eq:2.1.7}) are finite.  The difference
between $C^{(k,l), \mu}(I)$ and $ C^{k+l, \mu}(I)$ norms lies in that
in the former we only take up to $k$-th derivatives for the variable
$t$, thus distinguishing the `time' dimension $t$ from the `spatial'
dimensions on $M$. We define the following space of connections,
\begin{equation}
\label{eq:2.1.8}
\cS^{k, \mu}(I) = \{ A + \be dt: A(t) \in C^{k, \mu}(I) , \be(t)
\in C^{(k-1, 1), \mu}(I) \}.
\end{equation}
and let the norm of $\cS^{k, \mu}(I)$ be
\begin{equation}
\label{eq:2.1.9} 
|A + \be dt |_{\cS^{k, \mu}(I)}  = |A|_{C^{k, \mu}(I)} +
 |\be|_{C^{(k-1, 1), \mu}(I)} 
\end{equation}
The reason for this definition will become clear in Section 2.2.

The natural space of gauges for connections in $C^{k, \mu}(M)$ is the
space $C^{k+1, \mu}(M)$ of gauges. The natural space of gauges for
connections in $\cS^{k, \mu}(I)$ is given by gauges in $C^{(k, 1),
\mu}(I)$, because we see from the formula
\begin{equation}
\label{eq:2.1.10}
g(A + \be dt) = g A g^{-1} - dg \, g^{-1} + (g\be g^{-1} - \frac{\partial g
}{\partial t} g^{-1}) dt
\end{equation}
that $C^{(k, 1), \mu}(I)$ gauges is the minimal space of gauges that
act on $\cS^{k, \mu}(I)$ continuously.

For notations of $L^2$ and Sobolev norms, by fixing smooth covariant
derivatives on $E$ and  on $E \times I$, we define
\begin{alignat*}{2}
\| A \| &= \left(\int_{M} |A|^2 d\s \right)^ {\frac 12}, &\quad
\| A \|_{H^{l}}  & = \left( \sum_{0 \leq i \leq l} \int_M |\nabla^i
A|^2 d\s \right)^{\frac12},\\
\| A(t) \|_I &= \left( \int_{M \times I} | A(t) |^2 dt d\s \right)^
{ \frac 12}, &\quad \| A(t) \|_{H^{l}(I)} &= \left( \sum_{0 \leq i+j
\leq l} \int_{M \times I} |\nabla^i \pt_t^j A(t)|^2 dt d\s \right)^{\frac12}.
\end{alignat*}

Sometimes we abuse notation and use a connection $A$ on $M$ to
represent the time independent connection $p^*(A)$ on $M \times I$;
the meaning should be clear from the context. Gauge transformations
are usually assumed to be in the natural space of gauges for the
connections they act on; hence by $\Gamma(\Aut P)$ we often mean
$C^{k+1, \mu}(M, \Aut P)$, etc. If $A$ is a connection on $M$, we let
$\Stab(A) = \{ \s \in \Gamma(\Aut P) : \s (A) = A \} $ be the
stabilizer of $A$ in the gauge group.

Because the cylindrical coordinates $\phi$ is a conformal map, the
quadratic growth hypothesis on the curvature in
Theorem~\ref{th:1.1} translates to the uniform bound condition on
curvature
\begin{equation}
\label{eq:2.1.11}
|F_{\tilde{A}}| = (|F_A|^2 + |\dot{A} - d_A \be|^2)^{1/2} \leq C,
\end{equation}
for a connection $\tilde{A}$ on $E \times [0, \infty)$.  For a
Yang-Mills connection $\tilde{A}$ with \eqref{eq:2.1.10},
Prop~\ref{prop:1.3.2} implies the following compactness property,
\begin{lemma}
\label{lem:2.2}
Given $L>0$. For any sequence of positive numbers $\{R_i\}$ going to
infinity, there exists a subsequence ${i'} $, gauge transformations
$\s_{i'} \in \Gamma(\Aut P\times[R_{i'}, R_{i'}+L])$ and a Yang-Mills
connection $A$ on $M$ of $\tilde{A}$, such that
$$ \lim_{i' \ra \infty} |\s_{i'} (\tilde{A}) - A|_{C^{k,\mu}([R_{i'}, R_{i'} + L])} = 0 $$
\end{lemma}
%

\section{Constructions of Gauges}
\label{sec:2}
In this section we constructed the `standard form' gauges for Yang-Mills
connections with an isolated singularity in cylindrical coordinates
with certain estimates on norms of the connection. These estimates
will be important for us to prove the long-time existence of this gauge and
convergence of the connection in the next section. 
\subsection{Some bounds of gauges}
\label{sec:2.2} 
We give two well-known bounds about norms of gauges and connections in
the following lemma.
\begin{lemma}
\label{lem:2.3}
Let $\tilde{A}_0$ and $\tilde{A}_1$ be connections in $C^{k,
\mu}(I)$. $k \geq 0$, $0 \leq \mu \leq 1$. Assume that $\s$ is a gauge
on the bundle $E \times I$, we have\\
a) If $\max \{
|\tilde{A}_0|_{C^{k,\mu}(I)}, |\s (\tilde{A}_0)|_{C^{k, \mu}(I)} \}
\leq C$, then there exists $C_1 = C_1(k, C)$ such that
$$ |\s|_{C^{k+1, \mu}(I)} \leq C_1$$
b) There exists $C = C(k)> 0$ independent of $A_i$ and $\s$ such that
$$ |\s(\tilde{A}_1) - \s(\tilde{A}_0)|_{C^{k, \mu}(I)} \leq  C
(|\s|_{C^{k, \mu}(I)})^2 |\tilde{A}_0 - \tilde{A}_1|_{C^{k, \mu}(I)} $$
\end{lemma}
\begin{proof}
a) We have
\begin{equation}
\label{eq:2.4'}
d_{\tilde{A}_0}\s = (\s(\tilde{A}_0) - \tilde{A}_0)\cdot \s
\end{equation}
We also observe the fact that $|\s|_{C^0(I)} \leq C$ because the
structure group $G$ of $\Aut P$ is compact. Therefore we can apply a
bootstrapping procedure to (\ref{eq:2.4'}) and prove the required
estimates. This estimate is standard, see for example \cite[2.3.7]{DK}
.\\ b) We have
\begin{equation}
\label{eq:2.2.0} 
\s(\tilde{A}_1) - \s(\tilde{A}_0) = \s \cdot ( \tilde{A}_1 - \tilde{A}_0 ) \cdot \s^{-1}
\end{equation}
Notice that $d(\s^{-1}) = -\s^{-1}\,d\s\,\s^{-1}$ and therefore,
\begin{equation}
\label{eq:2.2.0'}
|\s^{-1}|_{C^{k, \mu}(I)} \leq C(k) |\s|_{C^{k, \mu}(I)}
\end{equation}
Combining (\ref{eq:2.2.0}) and (\ref{eq:2.2.0'}), we have the required bounds.
\end{proof}
\begin{remark}
The above lemma also holds if we look at other suitable spaces of
connections and natural spaces of gauges acting on them. For example,
the conclusion of the lemma hold if we replace the $C^{k, \mu}(I)$
norm for connections and the $C^{k+1, \mu}(I)$ and $C^{k, \mu}(I)$
norms for gauges by the $\cS^{k, \mu}(I)$ norm for connections and the
$C^{(k, 1), \mu}(I)$ and $C^{(k,0), \mu}$ norms for gauges
respectively. Similarly, the lemma holds for the $C^{k,\mu}(M)$ norm
for connections and the $C^{k+1,\mu}(M)$ and $C^{k, \mu}(M)$ norm of
gauges. We shall frequently use these bounds in the following, often
implicitly.
\end{remark}
Now we assume that $\tilde{A}$ is a smooth Yang-Mills connection on $E
\times [t_0, \infty)$ with $ \sup | F_{\tilde{A}}| \leq C$. In the
following lemma we give a bound on the set of tangent connections of
$\tilde{A}$. Here we regard tangent connections of $\tilde{A}$ as
Yang-Mills connections on $E$ over $M$ by time independence.
\begin{lemma}
\label{lem:2.3a}
Let $\cC$ = \{the set of tangent connections of $\tilde{A}$ \}, then there
exists $c_1 = c_1(k, \tilde{A}) > 0$ such that for any $A \in \cC$,
there exists $g \in C^{k+1, \mu}(\Aut P)$ with $|g(A)|_{C^{k, \mu}} \leq c_1$.
\end{lemma}
\begin{proof}
If $A_0 \in \Gamma(M, \wedge^1(\fg_E))$ is a tangent Yang-Mills
connection of $\tilde{A}$, then the bound on the curvature of
$\tilde{A}$ and the smooth convergence in Prop~\ref{prop:1.3.2} imply
that $\sup_{x \in M} |F_{A_0}(x)| \leq C$.  Let $A_i$ be a sequence of
tangent cones for $\tilde{A}$, then we may apply Prop~\ref{prop:1.3.2}
to obtain a subsequence $A_{i'}$ and gauge transformations $g_{i'}$,
such that $g_{i'}(A_{i'})$ converges smoothly to a tangent cone $A_0$
of $\tilde{A}$. This compactness property implies the bound in the
lemma, for any $k \geq 0$.
\end{proof}

The following lemma will be useful in Section 2.3.
\begin{lemma}
\label{lem:2.3b}
Given $L>0$, $\e>0$, there exists $R_0 = R_0(\e, L) >0$ such that if $R
\geq R_0$, then there exists a gauge transformation $g \in \Gamma(\Aut P
\times [R, R+L])$ and a tangent connection $A$ of $\tilde{A}$, such that
\begin{align}
\label{eq:2.4.7}
&|g(\tilde{A}) - A|_{\cS^{k, \mu}([R, R+L])}  \leq \e,\\
\label{eq:2.4.8}
&|A|_{C^{k, \mu}}\leq c_1, 
\end{align}
where $c_1 = c_1(k, \tilde{A})$ is the constant given in
Lemma~\ref{lem:2.3a}.
\end{lemma}
\begin{proof}
Assume that the lemma is not true, then there exist $\e_j \ra 0$ and $R_j
\ra \infty$ such that there does not exist tangent Yang-Mills
connection $A$ such that (\ref{eq:2.4.7}) and (\ref{eq:2.4.8}) are true
for $R$ replaced by $R_j$. By Lemma~\ref{lem:2.2}, there
exist a subsequence $\{ j' \}$ and $g_{j'} \in \Gamma(\Aut P \times
[R_{j'}, R_{j'}+L])$ such that $g_{j'}(\tilde{A})|_{[R_{j'}, R_{j'}+L]}
\stackrel{\cS^{k, \mu}}{\ra} A$, a tangent connection of $\tilde{A}$. By virtue
of Lemma~\ref{lem:2.3a}, up to a further $C^{k+1, \mu}$ gauge
transformation on $M$, we may assume $ |A|_{C^{k, \mu}} \leq c_1$. Now taking
$j'$ large will give us a contradiction.
\end{proof}
%
%
\subsection{Connections in standard form}
\label{sec:2.3} 
For Yang-Mills connections, there is a standard way of fixing gauges,
i.e. the {\em Coulomb gauge} (also called the {\em Hodge gauge}). We
say that $B$ is in the Coulomb gauge relative to $A$ if
$$d_A^*(B-A) = 0$$ (see \cite[2.3.1]{DK} or \cite{FU} for
example). The above Coulomb gauge equation and the Yang-Mills equation
(\ref{eq:1.1}) form an elliptic system. In the following, we find a
suitable gauge which gives rise to an elliptic system for a connection
$\tilde{A}$ on $E \times I$, $I$ being an interval of possibly infinite
length. Our choices of gauges are based on those in Morgan, Mrowka
and Ruberman \cite[2.4.3]{MMR}. The norms in the following proposition are
natural to our setting. 
\begin{proposition}
\label{prop:2.3.1}
Assume that $A_0$ is a smooth connection on $E$ and $\tilde{A}= A(t) +
\be(t) dt $ is a smooth connection on $E \times I$ . There exist $\e_1
= \e_1(A_0) >0$ and $C = C(A_0) >0$ such that if
\begin{equation}
\label{eq:2.3.0}
|\tilde{A} - A_0|_{\cS^{k, \mu}(I)} =  |A(t) - A_0|_{C^{k,
 \mu}(I)} + | \be(t) |_{C^{(k-1, 1), \mu}(I)} \leq \e_1,
\end{equation}
then there exists a gauge transformation $g \in C^{(k,1), \mu}(I)$,
such that $g (\tilde{A}) = A_1(t) + \be_1(t)dt$
satisfies
\begin{align}
\label{eq:2.3.1}
& d^*_{A_0}(A_1(t) -A_0) = 0,\\
\label{eq:2.3.2}
& \be_1(t) \in  \Ker (d_{A_0})^{\perp},\\
\label{eq:2.3.3}
&|g(\tilde{A}) - A_0|_{\cS^{k, \mu}(I)} = |A_1(t) - A_0|_{C^{k,
 \mu}(I)} + | \be_1(t) |_{C^{(k-1, 1), \mu}(I)} \leq C(1+|I|)^2
 |\tilde{A} - A_0|_{\cS^{k, \mu}(I)},\\
\label{eq:2.3.3a}
&|g - \Id|_{C^{(k,1),\mu}(I)} \leq C (1
 +|I|) |\tilde{A} - A_0|_{\cS^{k, \mu}(I)}. 
\end{align}
where  $\Ker (d_{A_0})^{\perp} = \Ker
(d_{A_0}: \Gamma(\fg_E) \ra \Omega^1(\fg_E))^{\perp} \subset
\Gamma(\fg_E)$. Furthermore, $g$ is unique up to the pullback of an
element of $\Stab(A_0)$, i.e., up to the composite with $p^*(\tau) \in
\Gamma(\Aut P \times I)$, where $p: \Aut P \times I \ra \Aut P$ is the
projection and $\tau \in \Stab(A_0) \subset \Gamma(\Aut P)$.
\end{proposition}

Following the terminology in \cite{MMR}, we make the following
 definition
\begin{definition}
We call the choice of gauge
$g(\tilde{A}) = A_1(t) + \be_1(t) dt$ satisfying \eqref{eq:2.3.1} and
\eqref{eq:2.3.2} in the above lemma a {\em standard form}  of
$\tilde{A}$ around $A_0$.
\end{definition}
Notice that the condition $d^*_{A_0} (A(t)-A_0) =
0$ means $A(t)$ is in Coulomb gauge relative to $A_0$. In the lemma,
if $|I| = \infty$, then (\ref{eq:2.3.3}) is vacuous.

\begin{proof}
We shall find the desired gauge in two steps.  \\
{\em Step 1.\,} We find a
gauge under which the Coulomb gauge condition (\ref{eq:2.3.1}) and the
bounds \eqref{eq:2.3.3} and \eqref{eq:2.3.3a} are satisfied. Consider
the gauge action map at $A_0$,
\begin{align*}
\Phi:\, C^{k+1, \mu}(\Aut P) &\times C^{k, \mu}(\Lambda^1(
\fg_E)) \ra  C^{k, \mu}(\Lambda^1( \fg_E))\\
&(g, ~~a)  \mapsto  - d_{A_0}g\, g^{-1}+  g\,a\,g^{-1}  = g(A_0 +a) - A_0
\end{align*}
The differential of this map at $(0, 0)$ is
\begin{align*} 
D\Phi:  C^{k+1, \mu}(\fg_E) & \times C^{k, \mu}(\Lambda^1 (
\fg_E)) \ra  C^{k, \mu}(\Lambda^1 (\fg_E))\\
& (h, ~~b)  \mapsto  - d_{A_0}h + b
\end{align*}
We can write down explicitly a right inverse of $D\Phi$ 
$$(D\Phi)^{-1}: C^{k, \mu}(\Lambda^1 (\fg_E)) \ra (C^{k+1, \mu}(\fg_E)
\cap \Ker(d_{A_0})^{\perp}) \times (C^{k, \mu}(\Lambda^1 (\fg_E)) \cap\Ker(d_{A_0}^*) )$$
by $(D\phi)^{-1}(a) = (h, b) $, where
$$h = (d_{A_0}^*d_{A_0})^{-1}(d_{A_0}^*a), \quad b = a + d_{A_0}h. $$
Here we used the fact that
\begin{equation}
\label{eq:2.3.5}
d_{A_0}^* d_{A_0}: C^{k+1, \mu}(\fg_E)\cap \Ker(d_{A_0})^{\perp}
\ra C^{k-1, \mu}(\Lambda^1 (\fg_E))\cap \Im (d_{A_0}^*)
\end{equation}
is invertible and has a continuous inverse.  Using the implicit
function theorem, we see that a right inverse $\Psi$ of $\Phi$ can be
defined in a neighborhood of $(0,0)$. It is easy to see that map
$$ \Psi:\, C^{k, \mu}(\Lambda^1(\fg_E)) \supset U \ra C^{k+1,
\mu}(\Aut P) \times C^{k, \mu}(\Lambda^1( \fg_E) \cap
\Ker(d^*_{A_0})) $$ is of class $C^k$ (note that the space of
connections is affine here). If $\Psi(a) = (g, b)$, then $A_0 + a =
g(A_0 + b)$ and $b \in \Ker(d^*_{A_0})$. This $\Psi$ gives the
well-known local slice structure for the configuration space of
connections under gauge actions.

Now we assume that $\tilde{A} = A(t) + \be(t)dt = A_0 + a(t) + \be(t)dt$
satisfies (\ref{eq:2.3.0}). In particular, $|a(t)|_{C^{k, \mu}}
\leq \e_1$, $\forall t \in I$. If $\e_1$ is sufficiently small, for
each fix $t$, we may define gauge transformation $g(t)$ by requiring
that $g(t)(A_0 + a(t)) = A_1(t) = A_0 + a_1(t)$ with $a_1(t) \in
\Ker(d^*_{A_0})$, i.e., by letting
$$ \Psi(a(t)) = (g^{-1}(t), a_1(t)).$$
Although $g(t)$ is obtained
separately for each $t$, we may view $g$ as a gauge transformation on
$E \times I$, and we have 
$$ g(\tilde{A}) = A_0 + a_1(t) + (g(t) \be(t) g(t)^{-1}- \frac{\pt}{\pt t}
g(t) g^{-1}) dt $$
Now that (\ref{eq:2.3.1}) is satsified, we need to
show that $g \in C^{(k, 1), \mu}(I) $ and the bounds (\ref{eq:2.3.3})
and (\ref{eq:2.3.3a}) hold. We observe first that by the continuity of
$\Psi$ on H\"older spaces,
\begin{equation}
\label{eq:2.3.5.0.1}
\sup_{t \in I} |g(t) - \Id|_{C^{k+1), \mu}} \leq C \sup_{t \in I}|
 a(t)|_{C^{k, \mu}} \leq C |\tilde{A} - A_0|_{\cS^{k, \mu}(I)}.
\end{equation}
By differentiating the identity $\Psi(a(t)) = (g^{-1}(t), a_1(t))$,
we have 
\begin{equation}
\label{eq:2.3.5.1}
-g ^{-1} \frac{\pt}{\pt t} g(t) g^{-1} = D_1\Psi_{a(t)}( \frac{\pt}{\pt t} a(t) )
\end{equation}
Since $ \frac{\pt}{\pt t} a(t) \in C^{k-1, \mu}(\Lambda^1 (
\fg_E))$ and
$$ 
D_1\Psi : C^{k-1, \mu}(\Lambda^1 (\fg_E)) \ra  C^{k,
\mu}(\Aut P)
$$
is a continuous map, we have 
\begin{equation}
\label{eq:2.3.5.2}
\sup_{t \in I} | \frac{\pt}{\pt t}g|_{C^{k, \mu}} \leq C \sup_{t \in
I}|\frac{\pt}{\pt t} a(t)|_{C^{k-1,\mu}} \leq C |\tilde{A}
-A_0|_{\cS^{k, \mu}(I)}
\end{equation}
By taking further derivatives of (\ref{eq:2.3.5.1}) and using
\eqref{eq:2.3.0}, we have 
\begin{equation}
\label{eq:2.3.5.3}
\sup_{t\in I} | \frac{\pt^j}{\pt t^j}g|_{C^{k+1-j, \mu}}  \leq C |\tilde{A}
  -A_0|_{\cS^{k, \mu}(I)}, \quad \hbox{for $1 \leq j \leq k$.} 
\end{equation}
(\ref{eq:2.3.5.0.1}) and (\ref{eq:2.3.5.3}) now imply \eqref{eq:2.3.3}
and \eqref{eq:2.3.3a}. This completes the first step.

{\em Step 2.\,} By the first step, we may assume $\tilde{A} = A(t) +
\be(t)dt = A_0 + a(t) + \be(t)dt$, $a(t) \in \Ker(d_{A_0}^*)$ and
(\ref{eq:2.3.0}) holds. We decompose $\be(t)$ as follows,
$$ \be(t) = \be_0(t) + \be_1(t),~ \be_0(t) \in \Ker(d_{A_0}),~ \be_1(t)
\in \Ker(d_{A_0})^{\perp}. $$
Since $\be(t) \in C^{(k-1, 1), \mu}$, by (\ref{eq:2.3.5}), we have 
$\be_1 = (d_{A_0}^*d_{A_0})^{-1}(d_{A_0}\be) \in C^{(k-1,1), \mu}(I)$
and $\be_0 = \be - \be_1 \in C^{(k-1, 1), \mu}(I)$ and
\begin{equation}
\label{eq:2.3.6a}
|\be_1(t)|_{ C^{(k-1, 1), \mu}(I)} + |\be_0(t)|_{ C^{(k-1, 1), \mu}(I)}
 \leq C |\tilde{A} - A_0|_{\cS^{k, \mu}(I)}.
\end{equation}
By standard ODE theory, the following linear
ordinary differential equation
\begin{equation}
\label{eq:2.3.6}
\left\{ \begin{array}{ll}
        &\frac{\pt }{\pt t}g = g \be_0(t), ~~\forall t \in I\\
        & g (0) = \Id.
        \end{array}
\right.
\end{equation}       
has a solution $g \in C^{k, \mu}(\Aut P \times I)$ (note $g(t)$
have one more derivative in $t$ than $\be_0(t)$ ).
Since $\be_0(t) \in \Ker(d_{A_0})= \rm{ the Lie algebra of }
\Stab(A_0)$, it follows that $g(t) \in \Stab(A_0), \forall t \in
I$. i.e. $g(t)(A_0) = A_0$, or equivalently,
\begin{equation}
\label{eq:2.3.6.1}
d_{A_0} g (t) = 0.
\end{equation}
Applying $g= g(t)$ to $\tilde{A}$, we have by (\ref{eq:2.3.6}), 
$$g(\tilde{A}) = g(t)(A(t)) + (g \be(t) g^{-1}- \frac{\pt g}{\pt t}
g^{-1})dt = A_0 + a_1(t) + g(t)(\be_1(t)) dt,$$
where $a_1(t) \in
\Ker(d_{A_0}^*)$ and $g(t)(\be_1(t)) = g(t) \be_1(t) g(t)^{-1}$. The
following lemma implies that $g (\be_1(t)) \in \Ker(d_{A_0})^{\perp}$,
hence $g(\tilde{A})$ is in the standard form around $A_0$.

\begin{lemma}
Assume that $A_0$ is a connection on $E$, $g \in \Stab(A_0)$ and $A \in
A_0 + \Ker(d_{A_0}^*)$, then $$g(A) \in A_0 + \Ker(d_{A_0}^*).$$ If
$\be \in \Ker(d_{A_0})^{\perp} \subset \Gamma(\fg_E)$, then
$$ g(\be) = g \be g^{-1} \in \Ker(d_{A_0})^{\perp}.$$
\end{lemma}
\begin{proof}
Since $g(A_0) = A_0$, we have
\begin{align*}
d_{A_0}^*(g(A)- A_0) & = d_{g(A_0)}^*(g(A) - g (A_0)) = g\circ
d_{A_0}^* \circ g^{-1} (g (A - A_0) g^{-1})\\
& = g (d_{A_0}^*(A-A_0)) g^{-1} = 0
\end{align*}
Hence $g(A) \in A_0 + \Ker(d_{A_0}^*)$. Similar to above, we can prove
$g \in \Stab(A_0)$ preserves $\Ker(d_{A_0}) \subset
\Gamma(\fg_E)$. Since gauge action preserves the metric on $E$, $g$
also preserves $\Ker(d_{A_0})^{\perp}$.
\end{proof}

Next we show that the gauge $g =g(t)$ from \eqref{eq:2.3.6}
satisfies \eqref{eq:2.3.3} and \eqref{eq:2.3.3a}. Differentiating 
(\ref{eq:2.3.6}) on $M$ and integrating in $t$, we have
\begin{align}
\label{eq:2.3.6'}
|g -\Id|_{C^{k, \mu}(I)}& \leq C (|\be_0|_{C^{(k-1,1), \mu}(I)} +
 \int_I |\be_0|_{C^{(k-1,1), \mu}(I)} dt )\\
&\leq C(1+ |I|) |\tilde{A} - A_0|_{\cS^{k, \mu}(I)}. \nonumber 
\end{align}

By differentiation of (\ref{eq:2.3.6.1}) with respect to $t$, we have
\begin{equation}
\label{eq:2.3.6.2}
d_{A_0} \frac{\pt g(t)}{\pt t} = 0.
\end{equation}
The equations (\ref{eq:2.3.6'}), (\ref{eq:2.3.6.1}) and (\ref{eq:2.3.6.2}) and
the smoothness of $A_0$ enable us to bootstrap on the derivatives of
$g$ on $M$. Hence we
have $ g \in C^{(k, 1), \mu}(I)$ and the following estimates.
\begin{equation}
\label{eq:2.3.6.3} 
|g - \Id |_{C^{(k, 1), \mu}(I)} \leq C |g - \Id|_{C^{k,
 \mu}(I)} \leq C (1 + |I|) |\tilde{A} - A_0|_{\cS^{k, \mu}(I)}
\end{equation}

Because $g(t) \in \Stab(A_0)$ by \eqref{eq:2.3.6.1},
\begin{align*}
|g(\tilde{A}) - A_0|_{\cS^{k, \mu}(I)} & = |g(t)(A(t)- A_0)|_{C^{k, \mu}(I)} +
|g\be_1g^{-1}|_{C^{(k-1,1), \mu}(I)}\\
\hbox{ by \eqref{eq:2.3.6a} } & \leq C (1 + |g -\Id|_{C^{k, \mu}(I)})^2
|\tilde{A} - A_0|_{\cS^{k, \mu}(I)} \\
\hbox{ by \eqref{eq:2.3.6'}} & \leq C(1+|I|)^2|\tilde{A} - A_0|_{\cS^{k, \mu}(I)} 
\end{align*}
The uniqueness of the standard form gauge (up to
$\Stab(A_0)$) is not hard and is left to the reader. 
\end{proof}

\begin{remark}
Note that in Step 1 of the proof above, we can not improve the number
of time derivatives of the gauge action $g$ to be more than that of
the connection, this is the reason why we work with $C^{(k,1),
\mu}(I)$ gauges and  $\cS^{k, \mu}$ connections, which at the
beginning may seem strange.
\end{remark}

\begin{lemma}
\label{lem:2.3.2}
Assume that $I$ is an interval of length $L$, $0< \e' \leq \t'\leq \tau
<1$, $\t$ is sufficiently small depending on $A_0$ and $L$,
$\tilde{A}$ is a connection on $E \times I$. $A_0, A_1$ are Yang-Mills
connections on $M$ with
\begin{align}
\label{eq:2.4.7.8} 
&|A_1 - A_0|_{C^{k,\mu}}  \leq \tau' \\
\label{eq:2.4.7.10} 
&|\tilde{A} - A_1 |_{\cS^{k, \mu}(I)}   \leq \e'
\end{align}
Then there exist gauges  $\tilde{\eta} \in C^{(k,1), \mu}(I \times \Aut P)$ and
$\eta \in C^{k+1, \mu}(\Aut P)$
such that $\tilde{\eta} (\tilde{A})$, $\eta{A_1}$ are respectively in standard
form and Coulomb gauge around $ A_0$, and
\begin{align*}
&|\eta(A_1) - A_0|_{C^{k,\mu}}  \leq C\tau'  \\
&|\tilde{\eta}(\tilde{A}) - \eta(A_1)|_{\cS^{k, \mu}(I)} \leq C\e'
\end{align*}
\end{lemma}
\begin{proof} 
We follow the two-step construction of standard form gauges around $A_0$ in
the proof of Prop.~\ref{prop:2.3.1} to obtain $\tilde{\eta}$. $\eta$ is
given by the usual Coulomb gauge around $A_0$, obtained in the first
step of the proof of Prop.~\ref{prop:2.3.1}. By keeping track of the
norms, it is easy to prove the lemma. We leave the details to the
reader.
\end{proof}

We state in the following lemma some elliptic estimates for
Yang-Mills connections in standard form gauges and Coulomb gauges.
\begin{lemma}
\label{lem:2.3.3}
a) Assume that $A_0$ is a smooth connection on $E$ over $M$ and
 $\tilde{A} = A(t) + \be(t)dt = A_0 + a(t) + \be(t)dt \in \cS^{k, \mu}(I)
 $ is a Yang-Mills connection in standard form around $A_0$. Then the
 following hold:\\
1) There exists $\e_1 = \e_1(A_0) > 0$, such that if
$|a(t)|_{C^{k, \mu}(I)} < \e_1 < 1 $, we have
\begin{equation}
\label{eq:2.3.7}
|\be(t)|_{C^{(k-1, 1), \mu}}(I) \leq C |a(t)|_{C^{k-1, \mu}(I)} |\dot{a}|_{C^{k-2, \mu}(I)}
\end{equation} 
for some constant $C = C(A_0, k) > 0$.
2) There exists
 $\e_2 = \e_2(A_0)> 0$ such that if $ |\tilde{A} - A_0|_{\cS^{k,
 \mu}(I)} < \e_2$, then for any $s \in (0, (b-a)/2)$, 
\begin{equation}
\label{eq:2.3.8}
|\tilde{A} - A_0|_{\cS^{k, \mu}(I_s)} \leq C(s) \| a(t) \|_{L^2(M
 \times I)},
\end{equation}
for some constant $C(s) = C(A_0, k, s) > 0$, where $I_s = [a+s, b-s]$
if $I = [a, b]$. \\
b) Assume that $A_1$ is a Yang-Mills connection on
$M$ in Coulomb gauge around a Yang-Mills connection $A_0$ then there
exists $\e = \e(A_0) > 0$, such that if $|A_1 - A_0|_{C^{1,\mu}} <
\e$, then
\begin{equation}
\label{eq:2.3.9} 
| A_1 - A_0 | _{C^{k,\mu}} \leq C(k) \|A_1 -A_0\|, \forall k \geq 0.
\end{equation}
\end{lemma}

\begin{proof}
{\em Proof of a). \, } Let $ \dot{a}(t) = \frac {\pt}{\pt t} a(t)
$. It follows from \eqref{eq:2.3.1} that $d_{A_0}^*\dot{a} = 0$.  From
one of the Yang-Mills equations (\ref{eq:2.2}), we have $d_A^*(\dot{a}-
 d_A\be) = 0 $, hence
\begin{equation}
\label{eq:2.3.10}
\Delta_A\be = d_A^*d_A \be = d_A^* \dot{a} = d^*_{A_0} \dot{a} -*[a,
*\dot{a}]
= - *[a,*\dot{a}]
\end{equation}

For fixed $t$, if $|A (t) - A_0|_{C^{k,\mu}} $ is sufficiently
small, then $\Delta_{A(t)} : C^{k, \mu} \cap \Ker (d_{A_0})^{\perp}
\ra C^{k-2, \mu} \Im (d_A^*)$ is invertible. We denote the inverse by
$G_A$. Hence
\begin{equation}
\label{eq:2.3.11} 
\be = - G_A (*[a, *\dot{a}]),
\end{equation}
Hence if $|a|_{C^{k, \mu}(I)} < e_1 < 1 $, we have from
(\ref{eq:2.3.11}) the estimates \eqref{eq:2.3.7}. We note that $\be(t)$
has one less derivative in $t$ than $a(t)$. We also note that $\be$ is of
quadratic nature in terms of $a$.

Substitute (\ref{eq:2.3.11}) into the Yang-Mills equation
(\ref{eq:2.1}) and use the Coulomb gauge condition (\ref{eq:2.3.1}), we have
\begin{align}
\ddot{a} - & (n-4) \dot{a} - d_A^*F_{A_0+a} - d_{A_0}d_{A_0}^* a
+ d_A( (G_A (*[a, *\dot{a}]))^{\cdot} \nonumber\\
\label{eq:2.3.12}
&-(n-4) d_AG_A(*[a, *\dot{a}]) ) + (-1)^{n+1} *[G_A(*[a, *\dot{a}]) , * d_A
G_A(*[a, *\dot{a}])] = 0.
\end{align}

We observe that the left hand side of (\ref{eq:2.3.12}) is an
(pseudo-differential) elliptic
operator of $a$ on $M \times I$ if $|a|_{C^{k, \mu}(I)} $ is
sufficiently small. If $|\tilde{A} -A_0|_{\cS^{k,\mu}(I)}$ is
sufficiently small, we have the following a priori interior estimates
for connections in standard form,
\begin{equation}
\label{eq:2.3.13}
|a(t)|_{\cS^{k, \mu}(I_s)} \leq C(s)\| a(t) \|_{I}, 
\end{equation}
 where $I_s =
[a+s, b-s]$ if $I= [a, b]$ and $s < \frac{b-a}{2}$. Combining the
estimate \eqref{eq:2.3.7}, we obtain the estimate \eqref{eq:2.3.8}.

{\em Proof of b). \, } The Yang-Mills equation and the Coulomb gauge
condition give
\begin{equation}
\label{eq:2.3.14}
d_{A_0+a}^* F_{A_0+a} + d_{A_0}d_{A_0}^* a = 0
\end{equation}
From which and bootstrapping, we obtain easily the estimate \eqref{eq:2.3.9}.
\end{proof}
%
%
\subsection[Standard form gauges with bounds]{Standard form gauges with bounds for Yang-Mills connections}
\label{sec:2.4} 
We shall prove the following key proposition, which provides us with a
good gauge (actually, the standard form around some tangent
connection) on an interval and a lower bound of the connection at the
end of the maximal existence interval of such a gauge. This design of
gauge is motivated by those in Cheeger and Tian \cite{CT}.

\begin{proposition}
\label{prop:4.1}
Assume $\tilde{A}$ is a Yang-Mills connection on $E
\times [t_0, \infty)$ with $\sup_{M \times [t_0, \infty)}|
F_{\tilde{A}}| \leq C$. We fix a tangent Yang-Mills connection $A_0$ of
$\tilde{A}$. Given $L>0$, $R_1>0$, there exists $\t_1=\t_1(A_0, L)> 0
$ and a function $\e_1 = \e_1(\t)$ satisfying $0 < \e_1(\t) < \t$ such
that for $0 < \tau \leq \tau_1$ and $0 < \e \leq \e_1(\tau)$, there
exist $R \geq R_1$, integer $2 < N \leq \infty$, and gauge
transformation $g \in \Gamma (\Aut P \times [R, R'])$, where $R' = R +
\frac{2}{3}NL$, such that the following (a) and (b) hold:\\ (a) $
g(\tilde{A}) = A(t) + \be (t) dt $ is of standard form around $A_0$. \\
(b) The following estimates hold,
\begin{align}
\label{eq:2.4.4}
&|g(\tilde{A})- A_0|_{\cS^{k, \mu}([R, R+2L])} = |A(t) - A_0|_{C^{k,
\mu}([R, R+2L])} + |\be(t)|_{C^{(k-1, 1), \mu}([R, R+2L])} \leq \e,\\
\label{eq:2.4.5}
&|g(\tilde{A})- A_0|_{\cS^{k, \mu}([R, R'])} = |A(t) - A_0|_{C^{k,
\mu}([R, R'])} + |\be(t)|_{C^{(k-1, 1), \mu}([R, R'])}  \leq \tau,\\
\label{eq:2.4.6}
& |\frac{\pt}{\pt t}g(\tilde{A})|_{\cS^{k-1, \mu}([R, R'])} =
|\frac{\pt }{\pt t} A(t)|_{C^{k-1,
\mu}([R, R'])} + |\frac{\pt}{\pt t} \be(t)|_{C^{(k-2, 1), \mu}([R, R'])}\leq \e,
\end{align}
In particular, $|A(t) - A_0|_{C^{k, \mu}([R, R'])} \leq \tau$ and $
|\frac{\pt}{\pt t} A(t) |_{C^{k-1, \mu}([R, R'])} \leq \e$. \\
\indent Furthermore, if $N$ is the maximal integer such that the above gauge
$g$ can be extended to $[R, R+\frac23 NL]$ and  (a) and (b) are satisfied,
and $N < \infty$, then we have \\
(c) 
\begin{equation}
\label{eq:2.4.6.0} 
\sup_{[R'-L, R']}|A(t) - A_0|_{C^{k, \mu}} \geq c_2
 \tau,
\end{equation}
for some constant $c_2 = c_2(A_0, L) \in (0, 1).$ 
\end{proposition}
\begin{remarks}
(1) The condition on $\e$ and $\t$ means that $\e$ is sufficiently small
relative to $\t$ and $\t$ is sufficiently small relative to $1$. As a
simplification, it can just be said as ``If $0 < \e << \t <<1$
depending on $A_0$, $L$, then \ldots''. We shall use this expression
later to simplify our statements.

(2) \eqref{eq:2.4.4} means that the connection is initially very close
to the tangent connection $A_0$. \eqref{eq:2.4.5} means that the
connection stays bounded to $A_0$ on the existence interval $[R,
R']$. \eqref{eq:2.4.6} means that the times derivative of the
connection $A$ keeps being very small on $[R, R']$. \eqref{eq:2.4.6.0}
means that the connection is at a distance away from $A_0$ at the end
of the existence interval of the constructed gauge. In other words, if
the connection stays within a certain distance to $A_0$, then the
gauge would exist up to infinity. As we shall see later,
(\ref{eq:2.4.6}) will play an important role in our proof of
Theorem~\ref{th:1.1} and is well worth our efforts here.
\end{remarks}
\begin{proof}
Assume given $L>0$ and $R_1 >0$.  We consider numbers
$\e_2, \e, \t$ with $0<\e_2< \e < \t$, which will be determined later.
Fix
$$R > \max\{ R_1, R_0( \e_2, 3L) + L\}, $$
with $R_0$ as in Lemma~\ref{lem:2.3b}. Choose intervals
$$I_i = [R+\frac{2iL}{3}-L, R+\frac{2iL}{3}], \quad 0
\leq i < \infty. $$
These intervals are chosen so that $|I_i| = L$, $|I_i \cap
I_{i+1}|= \frac 13 L$ and $I_i \cap I_{i+2} =
\emptyset$. We have the following claim.
\begin{claim}
\label{cl:2}
There exist constant $c_2=c_2(A_0, L)$, integer $2< K \leq \infty$,
Yang-Mills connections $A_i \in C^{k, \mu}$ in Coulomb gauge around
$A_0$ and gauge $g$ on $M \times \bigcup_{r=1}^{K}I_i$, such that
$g(\tilde{A}) = A(t) + \be(t)dt$ is in standard form around $A_0$,
\begin{align}
\label{eq:2.4.9}
&|g(\tilde{A}) - A_i|_{\cS^{k, \mu}(I_i)}\leq \e, ~~0 \leq i \leq K,\\
\label{eq:2.4.10}
&|A_i - A_0|_{C^{k, \mu}} \leq \frac{\tau}{2}, ~~0 \leq i \leq K,
\end{align}
and additionally, if $K$ is maximal as such and $K < \infty$, then
\begin{equation} 
\label{eq:2.4.11}
\sup_{t \in I_K} |A(t) - A_0|_{C^{k, \mu}} \geq c_2 \tau,
\end{equation}
\end{claim}

We first show that the claim implies our proposition. Assume the claim
is true, we take $N = K$. \eqref{eq:2.4.9} at $i=0$ gives
\eqref{eq:2.4.4}. Differentiation of \eqref{eq:2.4.9} with
respect to $t$ gives \eqref{eq:2.4.6}. Combining \eqref{eq:2.4.9} and
\eqref{eq:2.4.10} gives \eqref{eq:2.4.5}. \eqref{eq:2.4.11} gives
\eqref{eq:2.4.6.0}.

In the following, we prove the above claim by induction on $i$. For
$0\leq i \leq 3$, we note that by Lemma~\ref{lem:2.2}, there exists
gauge $g$ on $[R-L, R+2L] = \bigcup_{0\leq i \leq 3} I_i$ so that
$g(\tilde{A})$ is in standard form around $A_0$ and
\begin{align*} 
|g_0(\tilde{A}) - A_0|_{\cS^{k, \mu}([R-L, R+2L])}  &\leq  \e_2 \leq \e
\end{align*}
Hence we may let $A_i = A_0$, $0 \leq i \leq 3$. It is obvious that if
$\e \leq \frac{\t}{2}$, then (\ref{eq:2.4.9}), (\ref{eq:2.4.10}) are
satisfied for $0\leq i \leq 3$.

Assume that $i \geq 4$ and the claim is true for all integers $j$ such
that $j \leq i$. In the proof following, the constants $C$ we used
only depend on $L$ and $A_0$. Assume that the induction hypothesis
gives $$g(\tilde{A}) = A(t) + \be(t) dt$$ on $ [R-L, R+ \frac{2i}{3}L]$
and Yang-Mills connections $A_j$ for $0 \leq j \leq i$ on $M$ such
that \eqref{eq:2.4.9} and \eqref{eq:2.4.10} hold. If $\sup_{t \in I_i}
|A(t) - A_0|_{C^{k,\mu}} \geq c_2 \tau$, where $c_2= c_2(A_0, L)$ is to
be fixed later, then we are done by setting $K =i$. So we assume
\begin{equation}
\label{eq:2.4.11a} 
\sup_{t \in I_i} |A(t) - A_0|_{C^{k,\mu}} \leq c_2 \tau.
\end{equation}
We have by (\ref{eq:2.4.9}) and (\ref{eq:2.4.11a}) ,
\begin{equation}
\label{eq:2.4.7.-1}
|A_i - A_0|_{C^{k,\mu}} \leq 2 c_2 \tau,
\end{equation}
if $\e \leq c_2 \tau$.

Lemma~\ref{lem:2.3b} gives a tangent Yang-Mills
connection $A'_{i+1}$, and a gauge transformation $g'$
on $I_{i+1}$ such that
\begin{align}
\label{eq:2.4.7'}
&|g'{(\tilde{A})} - A'_{i+1}|_{\cS^{k, \mu}(I_{i+1})}  \leq \e_2,\\
\label{eq:2.4.8'}
&|A'_{i+1}|_{C^{k, \mu}}  \leq c_1, 
\end{align}
Let $t_i = R + \frac{2i}{3}L - \frac{L}{6}$ be the middle point of $I
= I_{i} \cap I_{i+1}$ and define $h = g(t_i) \cdot g'(t_i)^{-1}$. By
(\ref{eq:2.4.9}), (\ref{eq:2.4.10}), (\ref{eq:2.4.7'}) and
(\ref{eq:2.4.8'}), we have
\begin{align*}
&|g'(t_i)(A(t_i))|_{C^{k,\mu}} \leq
|A'_{i+1}|_{C^{k,\mu}} + C\e_2 \leq C,\\
&|h(g'(t_i)(A(t_i)))|_{C^{k,\mu}} = |g(t_i)(A(t_i))|_{C^{k,\mu}}
\leq |A_i|_{C^{k,\mu}} + C\e \leq C.
\end{align*}
Hence by Lemma~\ref{lem:2.3}, we have 
\begin{equation}
\label{eq:2.4.7.2} 
|h|_{C^{k+1, \mu}} \leq C
\end{equation}
Therefore,
\begin{align}
&|h(A_{i+1}) - A_i|_{C^{k, \mu}} \leq
|h(A_{i+1}) - h (g'(t_i) (A(t_i)))|_{C^{k,\mu}} \nonumber \\
\label{eq:2.4.7.3}
&\quad \quad+ | g(t_i)(A(t_i)) -
A_i |_{C^{k,\mu}} \leq C \e_2 + C' \e \leq C \e
\end{align}
Let $g'' = h g'$ on $I_{i+1}$, then by (\ref{eq:2.4.7.-1}) and
\eqref{eq:2.4.7'},
\begin{align}
\label{eq:2.4.7.5} 
&|g''(\tilde{A}) - h(A'_{i+1})|_{\cS^{k, \mu}(I_{i+1})} = |h \cdot
(g'(\tilde{A}) - A'_{i+1}) \cdot h^{-1}|_{\cS^{k,
\mu}(I_{i+1})} \leq C \e_2 \\
& |h(A'_{i+1}) - A_0|_{C^{k,\mu}}
\leq |h(A'_{i+1}) - A_i|_{C^{k,\mu}} + |A_i -
A_0|_{C^{k,\mu}}\nonumber \\
\label{eq:2.4.7.6}
&\qquad \leq C\e_2 + 2
c_2 \t \leq C c_2 \t 
\end{align}
if $\e \leq c_2 \t$.

We observe from (\ref{eq:2.4.7.5}) and (\ref{eq:2.4.7.6}) that $g''(\tilde{A})$
and $h(A'_{i+1})$ satisfies \eqref{eq:2.4.9} and \eqref{eq:2.4.10} in
the claim. But $g''(\tilde{A})$ is not in the standard
form around $A_0$ yet.  

In view of (\ref{eq:2.4.7.5}) and (\ref{eq:2.4.7.6}), we can apply
Lemma~\ref{lem:2.3.2} with $\tilde{A}$, $A_1$, $I$ and $\e'$, $\t'$ there
replaced by $g''(\tilde{A})$, $h A'_{i+1}$, $I_{i+1}$, $\e_2$ and $c_2
\tau$. Thus there exists gauges $\tilde{\eta}$ on $I_{i+1}$ and $\eta$
on $M$ respectively, such that
\begin{align*} 
&|\tilde{\eta} g''(\tilde{A}) - \eta (h (A'_{i+1})) |_{\cS^{k,
\mu}(I_{i+1})}  \leq c_3 \e_2\\
&|\eta h ((A'_{i+1})) - A_0|_{C^{k,\mu}} \leq c_4 c_2 \tau
\end{align*}
Compare $g$ and $\tilde{\eta} g''$ on $I = I_{i+1} \cap I_i$, we see
that they both make $\tilde{A}$ into standard form around $A_0$, hence
by the uniqueness of standard form gauges in Prop~\ref{prop:2.3.1},
they must differ by a pull back of $\sigma \in \Stab(A_0) \cap C^{k+1,
\mu}(M)$, i.e. $g=\sigma \tilde{\eta} g''$ on $I$. Therefore we can
extend $g$ to be over $I_{i+1}$ by letting
$$g|_{I_{i+1}} = \sigma \tilde{\eta} g''$$
Now we let $A_{i+1} = \sigma \eta h A'_{i+1}$, then  $g(\tilde{A})$
and $A_{i+1}$ are respectively in standard form and Coulomb gauge
around $A_0$, and
\begin{align*}
&|g(\tilde{A}) - A_{i+1}|_{\cS^{k, \mu}(I_{i+1})}  \leq C
|\tilde{\eta} g''(\tilde{A}) - \eta h A'_{i+1}  |_{\cS^{k, \mu}(I_{i+1})}   \leq c_5 \e_2 \\
& |A_{i+1} - A_0|_{C^{k,\mu}}  \leq C | \eta h A'_{i+1} -
A_0|_{C^{k,\mu}}  \leq c_6 c_2 \tau
\end{align*}
Finally, if we choose$c_2$, $\e_2$ sufficiently small such that $c_6 c_2
\leq \frac12$, $ c_5 \e_2 \leq \e$ and $0 << \e << \tau << 1$
satisfying all the requirements in the proof above, then
\eqref{eq:2.4.9} and \eqref{eq:2.4.10} are true for $i+1$ and 
the induction step of the claim follows. 
\end{proof}
%
%
\subsection{$L^2$ estimates of solutions at the end of maximal existence interval}
\label{sec:2.5} 
In the following proposition, we show that in the gauge constructed in
Prop.~\ref{prop:4.1}, near the end of the maximal existence interval
$[R, R']$ (if this is finite), the $L^2$ norms of the connection is
bounded away from zero and in average do not change much with respect
to time $t$.
\begin{proposition}
\label{prop:2.5.1} 
Given $L>0$ and $1>\eta >0$, if $0 < \e << \t <<1$ depending on $A_0$,
$L$ and $\eta$, and $\tilde{A} = A(t)
+\be(t)dt$ is a Yang-Mills connection and in the
gauge on $[R, R'] = [R, R+\frac23 NL]$ with respect to $\e, \t$ as in
Prop.~\ref{prop:4.1} and assume that $N < \infty$ is maximal there so that
(c) in Prop.~\ref{prop:4.1} holds, then
\begin{align} 
\label{eq:2.5.1} 
& \sup_{t\in [R'-L, R']} \|A(t) -A_0\| \leq (1 + \eta) \sup_{t\in
[R'-2L, R'-L]} \|A(t) -A_0\|,\\
\label{eq:2.5.2} 
& \sup_{t\in [R'-L, R']} \|A(t) -A_0\|  \geq c \t,
\end{align}
for some constant $c = c(L, \eta, A_0)$.
\end{proposition}
\begin{proof}
Since $N < \infty$, (c) in Prop.~\ref{prop:4.1} gives
\begin{equation}
\label{eq:2.5.1.5}
\sup_{t \in [R'-L , R']} |A(t) - A_0|_{C^{k,\mu}}  \geq c_2 \tau
\end{equation}
for $c_2 = c_2 (A_0, L)>0$. We first show the following claim.
\begin{claim}
\label{cl:2.5.1} 
Assume that`` $A(t) + \be(t) dt$ is a connection in the gauge of
Prop.~\ref{prop:4.1}, then there exists a Yang-Mills connection $A_1$ with
\begin{align}
\label{eq:2.5.1'} 
&\sup_{t \in [R'-2L , R']} |A(t) - A_1|_{C^{k,\mu}}  \le \d (\e)\\
\label{eq:2.5.1''}
&d_{A_0}^*(A_1 -A_0)=0
\end{align}
where $\d$ is an increasing function with $\lim_{\e \ra 0} \d(\e)=0$. 
\end{claim}
\begin{proof}[Proof of Claim]
Assume that the claim is not true, then there exists a sequence $\e_j \ra
0$ and a sequence of Yang-Mills connections $\tilde{A}_j = A_j(t) + \be_j(t)
dt$ on $I=[R'-2L, R'] $ in standard form gauge around $A_0$ and
satisfy conditions in Prop.~\ref{prop:4.1} (b), i.e.
\begin{eqnarray}
\label{eq:2.5.2'} 
&&|\tilde{A}_i - A_0|_{\cS^{k, \mu}(I)}  \leq \t,\\ 
\label{eq:2.5.3}  
&& |\frac{\pt }{\pt t}\tilde{A}_j|_{\cS^{k-1, \mu}(I)} \leq \e_j,
\end{eqnarray}
and there exists $\e_0 >0$ such that for any $j$ there does not exist a
Yang-Mills connection $A_1$ such that (\ref{eq:2.5.1'}) and
(\ref{eq:2.5.1''}) are true with $\delta(\e)$ and $A(t)$ replaced by $\e_0$ and
$A_j(t)$. By compactness of $\cS^{k, \mu}$ in  $\cS^k$, taking a
subsequence, we may assume $\tilde{A}_j \ra \tilde{A}' = A'(t) + \be'(t)dt$ in
$\cS^k(I)$ and it follows that $\tilde{A}'$ is also a Yang-Mills connection in
standard form around $A_0$ on $I$ and
$$ |\tilde{A}'|_{\cS^{k, \mu}(I)}  \leq \liminf_{i \ra \infty}
|\tilde{A}_j|_{\cS^{k, \mu}(I)}  \leq \tau $$

Taking limit as $j \ra \infty$ in (\ref{eq:2.5.3}), we have
$\dot{A'}(t) \equiv 0$, $\dot{\be}'(t) \equiv 0$. Hence, from the
Yang-Mills equation(\ref{eq:2.2}), we have
$$ d^*_{A'}d_{A'} \be' = 0 $$ Since $A'$ is close to $A_0$ in $C^k$
norm, $d^*_{A'}d_{A'}$ would be invertible on the space
$\Ker(d_{A_0})$ whileas $\be' \in \Ker(d_{A_0})$, hence $\be' \equiv
0$. From the other Yang-Mills equation (\ref{eq:2.1}), we have $A'$ is
Yang-Mills, hence $\tilde{A}'$ is in fact the pullback of a Yang-Mills
on $E$. Now (\ref{eq:2.1}) implies that
\begin{equation}
\label{eq:2.5.4}
d^*_{A_j(t)}F_{A_j(t)} = d_1(t) + d_2(t) + d_3(t),
\end{equation}
where
\begin{align*}
& \sup_{t \in I}|d_1(t)|_{C^{k-2,\mu}}  =  \sup_{t \in I}|\ddot{A_j}(t) +
\dot{A_j}(t)|_{C^{k-2,\mu}}
\leq 2 \e_j \hbox{ by (\ref{eq:2.5.3})}  \\
& \sup_{t \in I}|d_2(t)|_{C^{k-2,\mu}} =  \sup_{t \in I}|d_{A_j}\dot{\be}_j(t)|_{C^{k-2,\mu}} \leq
 C \e_j \hbox{ by (\ref{eq:2.5.3})} \\
& \sup_{t \in I}|d_3(t)|_{C^{k-2,\mu}} =  \sup_{t \in I}|(n-4)
d_{A_j}\be_j(t) + (-1)^{n+1} *[\be_j(t),
 *d_{A_j}\be_j(t)]|_{C^{ k-2,\mu}}  \\
& \quad \leq C  \sup_{t \in I} |\be_j(t)|_{C^{k-1,\mu}} := C \d_j 
\end{align*}
where $\d_j = \sup_{t \in I}|\be_j(t)|_{C^{k-1,\mu}} \ra 0$ since $\tilde{A}_j \ra A'$ in
$\cS^k(I)$. The left hand side of (\ref{eq:2.5.4}) is a uniformly elliptic
second order linear operator acting on $A_j(t) - A'$ with
coefficients bounded in $C^{k, \mu}(M)$ uniformly in $t$. Hence we can apply
elliptic estimates and bootstrapping to get
$$ \sup_{t \in I}|A_j(t) -A'|_{C^{k,\mu}} \leq C(\e_j + \d_j + \sup_{t
\in I} |A_j(t)
-A'|_{C^{k-1,\mu}} ) \ra 0 $$
But this is a contradiction to our assumption at beginning.
\end{proof}
Let $A_1$ be as in the above claim. Let $t_0 \in [R'-L, R']$
such that $|A_0 - A(t_0)|_{C^{k,\mu}} =
\sup_{t \in [R'-L, R']} |A_0 - A(t)|_{C^{k,\mu}} $, then by (\ref{eq:2.5.1.5}) and
(\ref{eq:2.5.1'}), we have
\begin{equation}
|A_1 -A_0|_{C^{k,\mu}}  \geq |A_0 - A(t_0)|_{C^{k,\mu}} - |A(t_0) -
A_1|_{C^{k,\mu}} \geq  c_2 \t - \d (\e ) \geq \frac12 c_2 \t
\end{equation}
if $\d(\e) \leq \frac12 c_2 \t$.
Now we have for any $ t \in [R'-2L, R']$,
$$
|A(t) - A_0|_{C^{k, \mu}} \geq |A_1 - A_0|_{C^{k, \mu}} - |A_1 -
A(t)|_{C^{k, \mu}} \geq \frac12 c_2 \tau -
 \d(\e) \geq \frac14 c_2 \tau
$$
 if $ \d(\e) \leq \frac{1}{4} c_2\tau$.
By a priori elliptic estimates for Yang-Mills connections in standard form
\eqref{eq:2.3.8}, we have
\begin{align*} 
\sup_{t \in [R'-2L, R'-L]}& \| A(t) - A_0 \| \geq C^{-1} |
A(R'-\frac{3L}{2}) - A_0|_{C^{k, \mu}} \geq \frac14 C^{-1}c_2\tau =
c_3 \tau\\
\sup_{t \in [R'-L, R']} & \| A(t) - A_0 \| \geq C^{-1} |
A(R'-\frac{L}{2}) - A_0|_{C^{k, \mu}} \geq c_3 \tau\\
\sup_{t \in [R'-L, R']} &\| A(t) -A_0 \| \leq \sup_{t \in [R'-2L,
R'-L]} \| A(t) - A_0 \|\\
\quad + 2CL|\frac{\pt}{\pt t} & A(t)|_{C^0([R'-2L, R'])} \leq (1+\eta) \sup_{t \in [R'-2L, R'-L]} \| A(t) - A_0 \| 
\end{align*}
if $4CL \e \leq \eta c_3\tau$.
This implies that if $0 < \e << \t <<1$ depending on $A_0$, $L$ and
$\eta$, then (\ref{eq:2.5.1}) and (\ref{eq:2.5.2}) hold.  
\end{proof}
%
%
\section[An asymptotic convergence result]{An asymptotic convergence
result for a class of evolution equations}
\label{sec:3}

By fixing the gauge as in the last section, we can reduce the main
theorem into an asymptotics problem for a certain nonlinear elliptic
evolution equation quite similar to the case treated in \cite[Theorem
1]{SL1}. The difference from the case there is that we only have
estimates and growth control (as given by \eqref{eq:2.4.6},
\eqref{eq:2.5.1} and \eqref{eq:2.5.2}) near the end of our existence
interval $[R, R']$, not on the whole interval.  However, to the up
side, we have the bound
$$|\frac{\pt^l}{\pt t} A(t)|_{C^{k-1, \mu}(I)} \leq \e$$ which allows
us to compare norms of connections at different points on the time
interval. We prove a general asymptotic convergence theorem relevant
to our case and complete the proof of Theorem~\ref{th:1.1} and
Theorem~\ref{th:1.2} in this section.

\subsection{A type of nonlinear evolution equations}
\label{sec:3.1}
Let $E$ be a vector bundle on
Riemannian manifold $M$ and let $\cE$ be a functional of `energy type'
defined for sections $a \in C^1(M, E)$ by
\begin{equation}
\label{eq:3.1.1}
\cE(a) = \int_M F(x, a, \nabla a)
\end{equation}
where $F = F(x, z, p)$ for $x\in M, \, z \in E$, $p
\in T_x M \otimes E_x$ depend smoothly on $(x, z, p)$ and $F$ is uniformly
convex in the $p$ variable for $p \in T_x M \otimes E_x$ and $|z|$,
$|p|$ small.  We also require that $F$ has analytic dependence on $(z,
p) \in E
\times T_x M \otimes E_x$ with uniform bounds
on $F$ and its derivatives in $z$, $p$ for sufficiently small $|z|$,
$|p|$ . By this we mean that there exists $c_0 >0$ such that
\begin{align*}
\label{eq:3.1.2}
F(x, z&+\sum_i \la_{1,i} w_i,\, p + \sum_j \la_{2,j} q_j) \\
&=
\sum_{|\alpha| \geq 0} F_\al (x, z, w_1, \ldots, w_m,
p, q_1, \ldots, q_{m'}) \la^\al,
\end{align*}
where $0 \leq i \leq m$, $0 \leq j \leq m'$, $m$ and $m'$ being the
dimension of fibers of bundle $E$ and $TM \otimes E$. $|z|, \, |w|, \,
|p_i|, \, |q_j|\leq c_0$, $z,\, w \in E_x$, $p_i,
\, q_j \in T_x M \otimes E_x$, $\la = (\la_{1,1}, \ldots, \la_{1,m},  \la_{2, 1}, \ldots, \la_{2, m'})
\in \RR^{m+m'}$, $|\la| \leq 1$. The above expansion and its
derivatives (in terms of $z$ and $p$ only) should converge absolutely
in the above domain of variables and there are uniform bounds
\begin{equation}
\label{eq:3.1.3} 
\sup_{|\al|+|\be| = j, \, |z|, \,|p| \leq c_0 } |D_z^{\alpha}D_p^{\beta} F (x, z, p)| \leq c(j).
\end{equation}
The Euler-Lagrange operator for $\cE(a)$, denoted by $\cM(a)$, is
uniquely characterized by
\begin{equation}
\label{eq:3.1.5}
-(\cM(a), b)_{L^2(M)} = \frac{d}{ds}\cE(a + s b)|_{s=0}.
\end{equation}
In other words, $\cM(a) = - \grad \cE(a)$). For simplicity, we also
require that
$$\cM(0) =0,$$
i.e. $0$ is a critical point of $\cE$. By
uniform convexity, we have that $\cM(a)$ is a second order
quasi-linear operator which is uniformly elliptic for $|a|_{C^1(M)}$
sufficiently small. Thus the linearization (the Jacobi operator at
$0$)
\begin{equation}
\label{eq:3.1.6}
La : = \frac{d}{ds}\cM(su)|_{s=0}.
\end{equation}
is a second order elliptic self-adjoint linear operator.

We shall fix $k \geq 2$ and $\mu > 0$ in this section for the
consideration of $C^{k, \mu}$ norms. Assume that $a (t) \in
C^{k, \mu}(M\times I, E\times I)$ is a section of $E \times I$. Fix $
\g > 0$ a constant and denote $\frac{\pt a(t)}{\pt t}$ by $\dot{a}(t)$. We consider the following equation
\begin{equation}
\label{eq:3.1.7}
\ddot{a}(t) - \g \dot{a}(t) + N(a(t)) + G_1(\dot{a}(t)) +
G_2(\ddot{a}(t)) =0.
\end{equation}
In the equation we require that $N(a)$ is a second-order quasi-linear
differential operator with smooth coefficients
and $N(a)$ is uniformly elliptic if $|a|_{C^1(M)}$ is sufficiently small and
$N$ approximates $\cM$ well in the following sense,
\begin{equation}
\label{eq:3.1.8} 
\| \cM(a) - N(a) \| \leq \frac14 \min \{ \|N(a)\|, \|\cM(a)\| \},
\end{equation}
for any $a \in C^{k,\mu}(M)$ with $| a |_{C^{k,\mu}} \leq \s$, where
$\s = \s(\cE)$ is a constant. Let $G_0(a): = N(a) - L(a)$ ($L$ as in
(\ref{eq:3.1.6})\,). We require in addition that
\begin{equation} 
\label{eq:3.1.9} 
G_i:\, C^{l, \mu}(M \times I) \ra  C^{l,
\mu}(M \times I) \quad \hbox{for }0 \leq i \leq 2,\,0\leq l \leq k-2
\end{equation}
are continuously differentiable linear maps and the derivatives
\begin{equation}
\label{eq:3.1.10} 
D^{\be}G_i :\;  C^{l, \mu}(M \times I) \ra C^{l, \mu}(M \times I),
\quad \forall  0\leq |\be| \leq k-2,\, 0 \leq l \leq
k-|\be|-2
\end{equation}
are continuous linear maps for $0\leq |\be| \leq k-2,\enskip 0 \leq l \leq
k-|\be|-2$ with the following uniform bounds on the operator norms:
\begin{equation}
\label{eq:3.1.11}
\| D^{\be}G_i(a)\|_{op} \leq C |a|_{C^{k,\mu}}, \quad \forall a \in C^{k, \mu}(M). 
\end{equation}
In fact, what we have in mind is that $G_i$'s are some pseudo-differential
operators which may have nonsmooth symbols.

If $a= a(t) \in C^{k, \mu}(M \times [t_1, t_2])$ is a solution of
(\ref{eq:3.1.7}) with $|a|_{C^{k, \mu}(I)} \leq C$, then it follows
from the standard Schauder theory and Sobolev theory (for example, in
\cite{CT} and \cite{Mo}), that for $\s \in (0, \frac12 (t_2 -t_1))$ and $0 \leq l
\leq k-2$,
\begin{align}
\label{eq:3.1.12}
&|a|_{C^{l,\mu}([t_1 + \s, t_2 -\s])} \leq c_l \s^{-(l + \mu)}\|
 a\|_{[t_1, t_2]}.
\end{align}
\noindent {\em Example.} If $\tilde{A} = A(t) + \be(t)dt = A_0 + a(t) +
\be(t)dt$ is a Yang-Mills connection in the standard form around a
Yang-Mills connection $A_0$. We claim that $a(t)$ satisfies an
equation of the form \eqref{eq:3.1.7}.

Let $\cE(a) = \YM(A_0 +a)$ in this case. We note that $\cE(a)$ is
uniformly convex if we restrict $a(t) \in \Ker(d_{A_0}^*)$
throughout. $\cE(a)$ is an analytic functional because the space of
connections is affine. We note that in contrast, the energy for
harmonic maps, which is analytic only when the metric on the target is
analytic. In the Yang-Mills case here, the linear fibers of the bundle
may be thought as the `target' space, which have
standard, hence analytic metrics. By
(\ref{eq:2.1}) and (\ref{eq:2.2}) and the fact $\be = - G_A (*[a,
*\dot{a}])$ ($G_A = (\Delta_A)^{-1}$, see the proof of
Lemma~\ref{lem:2.3.3}), $a(x, t)$ satisfies an equation of the form of
(\ref{eq:3.1.7}), where
\begin{align*}
\quad &\g = (n-4)\\
&N(a)= - (d_{A_0 +a}^* F_{A_0 +a} + d_{A_0}d_{A_0}^* a )\\
&L(a) = - ( d_{A_0}d_{A_0}^* a+ d_{A_0}^* d_{A_0} a+ (-1)^n *[a, *F_{A_0}]) \\
&G_0(a)= N(a) - L(a)= (-1)^n *[a, *(d_{A_0}a + a\wedge a)] + d^*_{A_0}(a \wedge a)\\
&G_1(\dot{a}) = d_A G_A (*[\dot{a}, *\dot{a}]) - (n-4)d_A G_A (*[a, *\dot{a}])\\
 \quad &+  (-1)^{n+1}*[ G_A (*[a, *\dot{a}]) , *d_A G_A (*[a, *\dot{a}])]\\
&G_2(\ddot{a})=  d_A G_A (*[a, *\ddot{a}])  
\end{align*}
We regard only one $\dot{a}$ in each term of $G_1$ as the variable,
others (dependent on $a$ and $\dot{a}$) will be seen as parts of
coefficients of $G_1$. Similarly, we only regard the $\ddot{a}$ in
$G_2$ as the variable. By the smoothing properties of the Green's
operator $G_A$, and (\ref{eq:3.3.10}), it is easy to see that the
operators $N, G_i$ satisfy the stated properties if $|a|_{C^{k,\mu}(I)}^* \leq
\t$ and hence the additional requirement following (\ref{eq:3.1.7})
are satisfied .

Assume that $\e$, $\t$, $\eta$, $L$ are positive constants, $\e \leq
\t$, $a= a(x, t) \in C^{k, \mu}(M \times [R, R'])$, where $R' -R \geq 3L$
and possibly $R'= \infty$.

\noindent {\it Definition.  }
We call $a = a(x, t)$  $(\e, \t, \eta, L)${\em-bounded}  on $[R, R']$ if
there exists $0< c = c(\eta, L)<1$ such that
\begin{align}
\label{eq:3.1.16} 
&|a|_{C^{k, \mu}([R, R'])} \leq \t \\
\label{eq:3.1.17} 
&|\frac{\partial a}{\partial t}|_{C^{k-1, \mu}([R, R'])} \leq
   \e \\
\label{eq:3.1.18} 
&|a|_{C^{k, \mu}([R, R+2L])} \leq \e
\end{align}
and additionally, if $R' < \infty$ then $a$ satisfies,
\begin{equation}
 \sup_{[R'-L, R']}\|a(t)\| \geq c \t.
\end{equation}

This definition, of course, is motivated by our estimates on
connections in Section 2.3 and Section 2.4.
We state our result about solutions to (\ref{eq:3.1.7}) as follows.
\begin{theorem}
\label{th:3.1.1} 
Fix $k \geq 5$ and $\mu > 0$. Given $L >0$, there exists $\eta \in
(0, 1)$ such that if $0 < \e << \t << 1$ depending on $E$, $L$ and
$\eta$, $a(t)$ is an $(\e, \t, \eta, L)$-bounded solution to
(\ref{eq:3.1.7}) and
\begin{equation}
\label{eq:3.1.21}
\cE(a(t)) - \cE(0) \geq - c \e, 
\end{equation}
then $R' = \infty$, $|a|_{C^{k, \mu}([R, \infty))}\leq \t$ and there
exists a $C^{2, \mu}$ critical point $\o$ of $\cE$ such that
\begin{equation}
\label{eq:3.1.22}
\lim_{t \ra \infty} |a(t) - w|_{C^{2,\mu}}
= 0
\end{equation}
\end{theorem}
%
%
\subsection{Growth estimates}
\label{sec:3.2}
Let $\cE$ be a functional as in the last subsection. Recall that $L$
is the linearization of the Euler-Lagrange operator $\cM = - \grad
\cE$. Assume that $\mu_1 \leq \mu_2 \leq \ldots$ and $\phi_1,\,
\phi_2, \dots$ are the complete set of eigenvalues and the corresponding
orthonormal (in $L^2$ norm) eigenfunctions of the operator $L$ on
$\Gamma(E)$. Let
\begin{equation} 
 \la_{i}^{\pm} = \frac{1}{2} (\g \pm \sqrt{\g^2 -
4\mu_i}). 
\end{equation}
Every solution of the linear evolution equation
\begin{equation}
\label{eq:3.1.12.5}
\cL (a) =  \ddot{a} - \g \dot{a} + L(a) =0,
\end{equation}
can be written as
\begin{align}
\label{eq:3.1.12.6} 
a (x,\, t) =& \sum_{i \in I_1} (a_i \cos \al_i t - b_i \sin
 \al_it)e^{\frac{\g t}{2}} \phi_i (x)\\
 & + \sum_{i \in I_2}(a_i + b_i t)e^{\frac{\g t}{2}} \phi_i(x) + \sum_{i \in I_3}(a_i
 e^{\la_i^{+}t}+ b_i e^{\la_i^{-}t}) \phi_i(x) \nonumber
\end{align}
for suitable constants $a_i$, $b_i$, where
\begin{align*}
I_1 &= \{ i: \mu_i < -\frac{\g}{4} \}, \quad \al_i = \Im \la_i^{+},\\
I_2 &= \{ i: \mu_i =  -\frac{\g}{4} \}, \\
I_3 &= \{ i: \mu_i >  -\frac{\g}{4} \}.
\end{align*}
And the $L^2$ norm square of $a(\cdot,\, t)$ can be written as 
\begin{align} 
\label{eq:3.1.14} 
\| a(t)\|^2 = &\sum_{i \in I_1}(a_i \cos \al_i t - b_i \sin \al_i t)^2
e^{\g t}\\
&\mbox{} + \sum_{i \in I_2}(a_i +b_i t)^2 e^{\g t} + \sum_{i \in I_3}(a_i
 e^{\la_i^{+}t}+ b_i e^{\la_i^{-}t})^2.  \nonumber
\end{align}
We let 
\begin{equation}
\label{eq:3.1.15} 
\d_1 =\min\{s: s\in\{\Re \la_i^{\pm}\}, s>0\}, \quad \d_2 =
\min\{|s|: s\in\{\Re \la_i^{\pm}\}, s<0\}
\end{equation}

We have
\begin{lemma}{\cite[Lemma 2]{SL1}}
\label{lem:3.2.1}
Assume $a \in C^{k, \mu}(M \times [0, 3L] )$ satisfy
(\ref{eq:3.1.7}), where $G_i$ satisfies properties above. For any given
$L>0, 1>\eta>0, \d < \frac14 \min\{ \d_1, \d_2 \}$, where $\d_1,\, \d_2$ as
in (\ref{eq:3.1.15}), there exists $\t_0 =
\t_0(L, \eta) >0$, such that if $|a|_{C^{k,\mu}([0, 3L])}  \leq \t_0$, then the following
are true. By denoting $S(j) = \sup_{t \in [(j-1)L,
jL]} \|a(t)\|$, we have
\begin{align*}
&\rm{(i) } S(2) \geq e^{\d L/2}S(1) \dra S(3) \geq e^{(\d_1 -\d)L}S(2)\\
&\rm{(ii) } S(2) \geq e^{-(\d_2 -\d)L}S(1) \dra S(3) \geq (1-\eta)S(2) \\
&\rm{(iii) } S(2) \geq e^{\d L/2}S(3) \dra S(1) \geq e^{(\d_2 -\d)L}S(2)\\
&\rm{(iv) } S(2) \geq e^{-(\d_1 -\d)L}S(3) \dra S(1) \geq (1-\eta)S(2) \\
&\rm{(v) If } S(2) \geq \max\{e^{-(\d_1 -\d)L}S(3), e^{-(\d_2 -\d)L}S(1)\},\\
&\rm{then } S(2) \leq (1+\eta) \inf_{t \in [L, 2L]} \|a(t)\|, \\
& \hbox{and }
\|\dot{a}(t) \|_{C^1} \leq \eta \| a(t)\|, \enskip \forall t\in [L, 2L].
\end{align*}
\end{lemma}
\begin{proof}[Sketch of proof]
First we can prove that (i)-(v) hold for solutions to the linear
equation (\ref{eq:3.1.12.5}) by using the expression of the $L^2$ norm
of (\ref{eq:3.1.14}). Then by a blow-up argument we can prove that for
$\t_0$ sufficiently small, the lemma holds for solutions to
(\ref{eq:3.1.7}). In the proof, we need the regularity properties of
solutions to (\ref{eq:3.1.7}).
\end{proof}

The next proposition is a simpler form of Theorem 4 in
\cite{SL1}.
\begin{proposition}
\label{prop:3.2.1}
Assume $a$ on $[R, R+NL]$ satisfies (\ref{eq:3.1.7}) as above. For
$0<\eta<1$, there exists $\t_0>0$, such that if
$0<|a|_{C^{k,\mu}([R, R+NL])}<\t_0$, then there exists integers $1
\leq k_1 \leq k_2 \leq N-1$ such that the following hold, where
$S(j)=\sup_{t \in [(j-1)L, jL]} \|a(t)\|$.\\
(a) $ S(j) \leq e^{-(\d_2 - \d)L} S(j-1)$, for $1\leq j
\leq k_1 -1$. \\
(b) $ \| a(t_1) \| \leq (1+\eta) \| a(t_2) \|$, for $t_1,\, t_2 \in
[k_1L, (k_2-1)L],\, \hbox{and }|t_1 - t_2| \leq L$ \\
\indent $ \| \dot{a}(t)\| \leq \eta \|a(t)\|$, for $ t \in [k_1L, (k_2-1)L]$.\\
(c) $S(j) \geq e^{(\d_1 - \d)L} S(j-1)$, for $k_2+1 \leq j
\leq N-1$.
\end{proposition}
\begin{proof}
We set
\begin{align*}
& k_1 = \min_{1\leq j\leq k-1} \{ S(j) \geq e^{-(\d_2-\d)L} S(j-1) \}
\\
& k_2 = \min_{k_1\leq j\leq N-1} \{ S(j+1) \geq e^{(\d_1-\d)L} S(j) \}
\end{align*}
And it is easy to check the theorem holds by repetitive use of the
previous lemma.
\end{proof}
This proposition allows us to conceptually divide the existence
interval into three parts according to the `growth rate' of the $L^2$
norms of $a$. In the first part which corresponds to case (a), the
$L^2$ norm of $a$ is decreasing exponentially (in average); in case
(b), the $L^2$ norm changes slowly in proportion; in case (c), it is
growing exponentially.

%
\subsection{Variational inequalities for analytical functionals}
\label{sec:3.3}
The inequalities in the following proposition are infinite dimensional
generalizations given by Leon Simon \cite{SL1} of {\L}ojasiewicz
inequalities with regard to critical points of analytic functions. Let
$\cE(a)$ be an analytic elliptic functional for $a\in C^{1}(M)$ as in
Section~3.1.
\begin{proposition}{\cite[Theorem 3]{SL1}}
\label{prop:3.3.1}
There exist constants $0<\theta<\frac12$, $2\leq \g$, $0 < \s$
depending only on $\cE$, such that if 
$|a|_{C^{k, \mu}} < \s$, then
\begin{equation}
\label{eq:3.3.1} 
\|\cM(a)\| \geq ( \inf_{\zeta \in \cS} \|a - \zeta \|)^{\g}
\end{equation}
where $\cS = \{ \zeta \in C^{k,\mu}(E): |\z|_{C^{k,\mu}(\Sigma)}< c_0,
\cM(\z) =0 \}$, and
\begin{equation}
\label{eq:3.3.2}  
\|\cM(a)\| \geq |\cE(a) - \cE(0)|^{1-\theta}
\end{equation}
\end{proposition}

Let $\s$, $\theta$ in this section be the same as in
Prop.~\ref{prop:3.3.1}. Now assume that $a \in C^{k,\mu} (M \times [t_1, t_2])$
satisfies the following equation
\begin{equation}
\label{eq:3.3.3} 
\dot{a} = N (a) + R(a)
\end{equation}
where $N (a)$ approximates well the gradient $\cM$ of $\cE$ in the
sense of (\ref{eq:3.1.8}) and
\begin{equation}
\label{eq:3.3.5} 
\| R (a)(t) \| \leq \frac 12 \| \dot{a}(t) \|, \quad \forall t \in [t_1,
t_2]
\end{equation}
We have 
\begin{lemma}{\cite[Lemma 1]{SL1}}
\label{lem:3.3.1}
Suppose $a$ satisfies $|a(t)|_{C^{k,\mu}([t_1, t_2])} \leq \s$, and
suppose that for some constant $\e > 0$,
\begin{equation}
\label{eq:3.3.6}
\cE(a(t)) > \cE(0) - \e ~~\hbox{for all } t \in [t_1, t_2]
\end{equation}
Then
\begin{equation}
\label{eq:3.3.7} 
\int_{t_1}^{t_2} \| \dot{a}(t) \| dt \leq C \theta^{-1} ( |\cE (a(t_1)) - \cE(0
)|^{\theta} + \e^{\theta}),
\end{equation}
In particular
\begin{equation}
\label{eq:3.3.8}
\sup_{t \in [t_1, t_2]} \| a(t) - a(t_1) \| \leq C \theta^{-1}( |\cE
(a(t_1)) - \cE( 0 )|^{\theta} + \e^{\theta}).
\end{equation}
\end{lemma}

The next lemma is a gauge invariant form of the inequality
(\ref{eq:3.3.2}) for connections and was essentially proven in \cite{MMR} .
\begin{lemma}
\label{lem:3.3.2}
Let $E$ be a vector bundle on $M$. $A$ is a  $C^{k, \mu}$
connection on $E$ and $B$ is a smooth Yang-Mills connection on $E$. There
exists $\e_3 > 0$ and $\theta \in (0, \frac 12)$, such that if $| A -B
|_{C^{k, \mu}} < \e_3$, then the following inequality holds,
\begin{equation}
\label{eq:3.3.9}
\left(\int_{M} |F_A|^2 - |F_B|^2 d\s \right)^{1-\theta} \leq 2 \| d_A^* F_A \|.
\end{equation}
\end{lemma}
\begin{proof}
If $\e_3$ is sufficiently small, by a gauge transformation to the
Coulomb gauge around $B$, we may assume that $A=B+a$ with $a \in
\Ker(d_B^*)$, and $|a|_{C^{k, \mu}} < C \e_3$. Let $U$ be a small
$C^{k, \mu}$ neighborhood of $0$ in $\Ker(d_B^*) \subset \Omega^1(\End
E|_{M})$ and consider $\cE: U \ra \RR$ by $\cE(b) = \YM(B+b)$,
$\forall b \in U$. Let $\cM $ be the Euler-Lagrange operator of $\cE$
on $U$. We claim that if $U$ is sufficiently small, $A = B+a$, $a \in
U$, then
\begin{equation}
\label{eq:3.3.10}
\| d_A^* F_A -\cM(a) \| \leq
\frac14 \min \{ \| \cM(a) \|, \, \| d_A^* F_A\| \}.
\end{equation}
Indeed, since $ \langle \cM(a), b \rangle_{L^2} = \langle d_{A}^*F_A,
\, b\rangle_{L^2}$, $\forall b \in \Ker(d_B^*)$ and $\cM(a) \in
\Ker(d_B^*)$, it follows that $\cM(a) = p_{\Ker(d_B^*)} (d_A^*F_A)$,
where $p$ is the $L^2$ projection. Let $G_B = d_B^* d_B ^{-1}$ as in
the proof of Lemma~\ref{lem:2.3.3}, we can easily see that
\begin{align} 
\label{eq:3.3.11}
\cM (a) = 
p_{\Ker(d_B^*)} (d_A^*F_A)& = d_A^*F_A - d_B G_B d_B^*(d_A^*F_A)\\
& = d_A^*F_A - d_B G_B (*[a, *d_A^*F_A])\nonumber,
\end{align}
where we used the identity
\begin{equation}
\label{eq:3.3.12}
d_A^* d_A^* F_A = \{ F_A , \, F_A \} =  0.
\end{equation}
where $\{\,,\,\}$ is defined by Lie bracket on the bundle parts and
Riemannian product on the form part and hence is skew-symmetric.  From
\eqref{eq:3.3.11} and the smoothing properties of $G_B$ it is easy to
derive (\ref{eq:3.3.10}) if $|a|_{C^{k,\mu}} $ is small.  Apply
(\ref{eq:3.3.2}) to $\cE$ and use (\ref{eq:3.3.10}) (note that $\cE$
is only defined on $\Ker(d_B^*)$, not on the space of all
sections; however Prop.~\ref{prop:3.3.1} is still true in this case
because $\Ker(d^*_B)$ is an analytic submanifold of the H\"older
spaces of sections), we obtain that there exists $\theta \in (0, \frac
12)$, such that
$$\| \cE(a) - \cE(0)\|^{1-\theta} \leq 2\| d_A^* F_A \|. $$
This last
inequality is exactly (\ref{eq:3.3.9}) which we want to prove.
\end{proof}

%
\subsection{Proof of Theorem~\ref{th:3.1.1}}
\label{sec:3.4}
Our method of proving Theorem~\ref{th:3.1.1} follows mostly the
methods of proving Theorem 1 in \cite{SL1} with modifications to our
case. We first apply the growth estimates in Section~\ref{sec:3.2} to
the time derivative $\dot{a}(t)$ of the solution to \eqref{eq:3.1.7}
and divide the interval into three parts. We estimate the integration
of $\| \dot{a}(t) \|$ on these three parts respectively and thus
obtain a bound on $\|a(t)\|$, which then gives long-term existence of
$a(t)$ and convergence. The estimates of $\|\dot{a}(t)\|$ on the
three parts use the estimates from Section~\ref{sec:3.3} and the
assumptions on the solution in Theorem~\ref{th:3.1.1}, especially the
condition \eqref{eq:3.1.17}.

Assume $R' < \infty$. We may assume that $R'=R +NL$ for $N \geq 2$ by
changing $R$ within an amount of $L$ if necessary. It is important to
note that by differentiation of (\ref{eq:2.1}), we have $\dot{a}$
satisfies an equation of (\ref{eq:3.1.7}) form on $[R, R']$, with the
same $\t$ (we need to change $C^{k, \mu}$ norms to $C^{k-1, \mu}$ norms for
$\dot{a}$; since $k\geq 5$, this is still in the regular range, an our
previous results hold without change). Note that $\ddot{a}$ in general
doesn't satisfy an equation in the form of \eqref{eq:3.1.7}.

Assume $\t < \t_0$, where $\t_0$ is as in Prop.~\ref{prop:3.2.1}. Then
by applying Prop.~\ref{prop:3.2.1} to $\dot{a}(t)$, we have $1\leq k_1
\leq k_2 \leq N-1$ for $\dot{a}$ such that the conclusions of
Prop.~\ref{prop:3.2.1} hold. We adopt the notation
\begin{equation}
\label{eq:3.4.1} 
S(j, \dot{a}) = \sup_{t\in [R+(j-1)L, R+jL]} \| \dot{a}(t) \|
\end{equation}
%

We shall assume in the rest of our proof that $\t \geq \e^\al$ for the
constant $\al = \frac{\th}{8}$, where $\th$ is the constant in Lemma
\ref{lem:3.3.1}.  We have the following claim,
\begin{claim}
\label{cl:3}
There exists a constant C depending on the functional $\cE$ such that
\begin{equation}
\label{eq:3.4.2}
\|a(t)\| \leq C \e^{2\al}, \quad \forall t \in [R, R+NL].
\end{equation}
\end{claim}
\begin{proof}[Proof of Claim]
First we consider the claim for $t \in [R, R+(k_1-1)L]$, by using
Prop.~\ref{prop:3.2.1} (a), we have
\begin{equation}
\label{eq:3.4.3} 
\| \dot{a}(t) \| \leq S([(t - R) /L]+1, \dot{a}) \leq 
e^{-(\d_2 -\d)(t-L)} S(1, \dot{a}) \leq e^{-\d(t-L)} \e
\end{equation}
Hence for $t \in [R, R+(k_1-1)L]$,
\begin{equation} 
\label{eq:3.4.4} 
\| a(t) \| \leq \|a(R)\| + \int_0^t \|\dot{a}\| dt \leq \e +
\d^{-1}e^{\d L}\e \leq \e^{\frac12}, 
\end{equation}
if $\e^{\frac12} \leq (1 +  \d^{-1}e^{\d L})^{-1}$.

Next we consider the case $t \in [R+(k_1-1)L, R+(k_2+1)L]$. By using
the fact that $|a(t)|_{C^{k, \mu}} \leq \t$ and $\| \ddot{a}(t) \|
\leq 1/8 \| \dot{a}(t)\|$, which follows from (b) in
Prop.~\ref{prop:3.2.1} if we take $\eta = 1/8$ there, we see that
Lemma \ref{lem:3.3.1} applies to $\dot{a}$ on $[R+k_1L, R+(k_2-1)L]$
as long as $\e$ and $\t$ are sufficiently small. Therefore
\begin{align}
\label{eq:3.4.5} 
\int_{R+k_1L}^{R+(k_2-1)L} \|\dot{a}\| dt & \leq C (|\cE(a(R+k_1L)) -
\cE(0) |^{\th} + \e^{\frac{\th}{2}} )\\
& \leq C (| a(R+k_1L) |_{C^1}^{\th} + \e^{\frac{\th}{2}}) \hbox{ by
\eqref{eq:3.1.1}}\nonumber \\
&\leq C (\sup_{t \in [R+k_1L -1, R+k_1L +1]} \| a(t) \|^{\th} +
\e^{\frac{\th}{2}}) \hbox{ by
\eqref{eq:3.1.12}}\nonumber \\
& \leq C \e^{\frac{\th}{2}} \leq C\e^{2\al} \hbox{ by
\eqref{eq:3.4.4}}\nonumber
\end{align}
and hence for $\t \in [R+(k_1-1)L, R+(k_2+1)L]$,  by
(\ref{eq:3.1.17})  $\|\dot{a}\| \leq C\e$ and we have
\begin{align}
\label{eq:3.4.6}
\| a(t)\| &\leq \| a(R+k_1L)\| + C \e^{2\al} + C \e\\
&\leq \|a(R+(k_1-1)L)\| + C\e + C\e^{2\al}\nonumber\\
&\leq C \e^{\frac12} + C\e + C\e^{2\al} \leq C\e^{2\al}. \nonumber
\end{align}
Finally, we consider the case $t \in [R+(k_2+1)L, R+NL]$.  It follows
from Prop.~\ref{prop:3.2.1} (c) that
$$ S(j, \dot{a}) \geq e^{(\d_1-\d)L}S(j-1, \dot{a}), \forall k_2 +1
\leq j \leq N-1. $$
We also know from \eqref{eq:3.1.17} that $S(N-1,
\dot{A}) \leq C\e$. Hence for any $t \in (R+(j-1)L, R+jL]$, where $k_2 +1
\leq j \leq N-1$, we have
$$\|\dot{a}(t) \| \leq S(j, \dot{a}) \leq e^{-(\d_1 - \d)(N-1-j)L}S(N-1,
\dot{a}) \leq Ce^{-(\d_1 - \d)(R + (N-2)L - t)}\e.$$
It follows that
\begin{equation}
\int_{R+(k_2+1)L}^{R+(N-1)L} \| \dot{a}(t)\|  \leq
\int_{R+(k_2+1)L}^{R+(N-1)L} Ce^{-(\d_1 - \d)(R + (N-2)L - t)}\e 
 \leq C \e.
\end{equation}
Hence for any $t \in [R+(k_2+1)L, R+NL]$, we have
\begin{equation}
\| a(t)\| \leq \| a(R+(k_2+1)L)\|  + C \e \leq C \e^{2 \al}
\end{equation}
\end{proof}

Now (\ref{eq:3.4.2}) implies in particular that
\begin{equation}
\label{eq:3.4.7} 
\sup_{t \in [R+(N-1)L, R+NL]} \| a (t) \| \leq  C\e^{2\al}
\end{equation}
which gives a contradiction to our assumption (\ref{eq:3.1.20}) if
$\e$ is sufficiently small relative to $\t$. Therefore we must have $R' =
\infty$. Hence by (\ref{eq:3.1.16}) , $|a|_{C^{2, \mu}([R, \infty))}
\leq \t$ and by \eqref{eq:3.1.17}, $\|\dot{a}(t)\| \leq \e$ for $ t
\in [R, \infty)$. It follows that $k_2 = \infty$, otherwise
Prop.~\ref{prop:3.2.1} (c) implies $\| \dot{a}(t)\|$ is going to
infinity for a sequence of $t$. Now \eqref{eq:3.4.4} and
\eqref{eq:3.4.5} imply that
$$\int_{R}^{\infty} \| \dot{a}(t) \| < \infty.$$ Hence there exits $w
\in L^2$ such that $a(t) \ra w$ in $L^2$ as $t \ra \infty$. Uniform
bounds on $| a(t) |_{C^{3, \mu}}$ implies that for a sequence $t_i \ra
\infty$, $a(t_i) \ra w$ in $C^{2, \mu}$. In fact, compactness implies
 that $a(t) \ra w$ in $C^{2, \mu}$ as $t \ra \infty$.  Hence $w
\in C^{2, \mu}$ and taking the limit of \eqref{eq:3.1.7} as $t_i \ra
\infty$, we see that $w$ is a critical point of $\cE$. This finishes
the proof of Theorem~\ref{th:3.1.1}.
%
%
\subsection{Proofs of Theorem~\ref{th:1.1} and Theorem~\ref{th:1.2}}
\label{sec:3.5} 
In this section we apply Theorem~\ref{th:3.1.1} to prove the convergence
of the connection to its tangent connection, hence the uniqueness of
tangent connections in Theorem~\ref{th:1.1}. We shall use the
monotonicity formula and Lemma~\ref{lem:3.3.2} of Yang-Mills
connections to show the desired rate of convergence in
Theorem~\ref{th:1.1}. The idea of using monotonicity formula and
Lemma~\ref{lem:3.3.2} comes from Leon Simon's work \cite[3.10 -
3.15]{SL2}, where energy minimizing harmonic maps with a tangent map
which has an isolated singularity is treated.

As before, assume that in the cylindrical coordinates, $\tilde{A} =
A(t) + \be(t)dt = A_0 + a(t) + \be(t)dt$, $t \in [R, R']$ is a
connection on $S^{n-1} \times [R, R']$, is in the standard form gauge
given in Prop.~\ref{prop:4.1}, then the conclusions of
Prop.~\ref{prop:4.1} and Prop.~\ref{prop:2.5.1} imply that $a(t)$ is a
$(\e, \t, \eta, L)$-bounded solution of an equation in the form of
\eqref{eq:3.1.7} on the interval $[R, R']$. The only condition in
Theorem~\ref{th:3.1.1} left is the energy lower bound
\eqref{eq:3.1.21}. This can be achieved by choosing the tangent
connection $A_0$ suitably at the beginning. We note that the set of
energies of tangent connections of $\tilde{A}$ is bounded from
below. We choose $A_0$ with energy very close to the infimum, then it
is easy to see \eqref{eq:3.1.21} must be satisfied. Hence
Theorem~\ref{th:3.1.1} applies to give the long-time existence of the
standard form gauge and the convergence of the connections $A(t)$ to a
tangent connection $A_0'$ as $t \ra \infty$. It is easy to see that
any other tangent connection must be gauge equivalent to this $A_0'$.

Next we proceed to show the rate of convergence for $A(t)+ \be(t)dt \ra
A_0$ as $t \ra \infty$. Without loss of generality, we may assume $R =
0$ and $\tilde{A}_1 =\phi^*(\tilde{A})$ is the original connection on
$B_1(0)\setminus \{ 0\}$. Since $\tilde{A}_1$ is stationary, by
monotonicity formula,
\begin{eqnarray}
\label{eq:3.5.-1} 
&&4\int_{B_{\rho}(0)} r^{4-n} | \frac {\pt}{\pt r} \rfloor F_{\tilde{A}_1} |^2
dx \nonumber \\
&\leq& \lim_{\s \ra 0} \{ \int_{B_{\rho}(0)}
\rho^{4-n} |F_{\tilde{A}_1} |^2 dx - \int_{B_{\s}(0)}
\s^{4-n} |F_{\tilde{A}_1} |^2 dx \} \nonumber\\
&=& \lim_{\la_i \ra 0} \{ \int_{B_{\rho}(0)}
\rho^{4-n} |F_{\tilde{A}_1} |^2 dx - \int_{B_{\la_i}(0)}
\la_i^{4-n} |F_{\tilde{A}_1} |^2 dx \}\nonumber\\
&=& \lim_{\la_i \ra 0} \{ \int_{B_{\rho}(0)}
\rho^{4-n} |F_{\tilde{A}_1} |^2 dx - \int_{B_1(0)}
 |F_{\tilde{A}_{1_{\la_i}}} |^2 dx \} \nonumber\\
&=&  \int_{B_{\rho}(0)}
\rho^{4-n} (|F_{\tilde{A}_1} |^2 - |F_{A_0}|^2) dx \nonumber \\
&\leq &  \frac {1}{n-4}  \int_{\pt B_{\rho}(0)}
\rho^{5-n} (|F_{\tilde{A}_1} |^2 - |F_{A_0}|^2 )d \s
\end{eqnarray}
where the last inequality follows from monotonicity formula
(\ref{eq:m1}) and the fact that $\frac{\pt}{\pt r} \lfloor F_{A_0} =
0$.  We note that there is no loss of curvature energy (on bounded
sets) by Prop.~\ref{prop:1.3.2}. Under our cylindrical coordinates,
let $\eta(t)=\dot{a}(t) - d_A \be(t)dt$ and $T = - \log( \rho) \geq R$,
by change of variables, (\ref{eq:3.5.-1}) becomes, for any $T \in [R,
R']$,
\begin{align}
\label{eq:3.5.0}
\int_T^{\infty} \| \eta(t) \|^2 dt &\leq
\frac{1}{n-4} \int_{S^{n-1}} | F_A(T) - \eta(T)dt|^2 - |F_{A_0}|^2 d\s\\
& = C \int_{S^{n-1}} | F_A(T) |^2 - |F_{A_0}|^2 d\s + C \| \eta (T)\|^2  \nonumber
\end{align}

Since $|A(t)-A_0|_{C^{k,\mu}} = |a(t)|_{C^{k,\mu}} \leq \t$, we can
apply Lemma~\ref{lem:3.3.2} to the right hand side of (\ref{eq:3.5.0})
and obtain
\begin{equation}
\label{eq:3.5.1}
\int_T^{\infty} \| \eta(t) \|^2 dt \leq  C\|d_{A(T)}^* F_{A(T)} \|^{\frac{1}{1-\theta}} + C \| \eta (T) \|^2, 
\end{equation}
From (\ref{eq:2.1}), we have
\begin{equation}
\label{eq:3.5.2}
\| d_{A(t)}^* F_{A(t)} \| \leq \| \ddot{a}(t) \| + \| \dot{a}(t) \|
+ \| d_A \dot{\be}(t) \| + \| d_A \be(t) \| + C \tau \| \be(t) \|
\end{equation}

\eqref{eq:2.3.8} implies that
\begin{equation}
\label{eq:3.5.8}
|\be(t)|_{C^{k, \mu}} + |\dot{\be}(t)|_{C^{k, \mu}} \leq C \t
|\dot{a}|_{C^{k, \mu}}.
\end{equation}
\eqref{eq:3.1.10}
The elliptic estimates \eqref{eq:3.1.12} implies that
\begin{equation}
\label{eq:3.5.9} 
|\ddot{a}(t)|_{C^{k,\mu}}+ |\dot{a}(t)|_{C^{k,\mu}} \leq C \| \dot{a} \|_{[t-1, t+1]} 
\end{equation}
Applying (\ref{eq:3.5.8}) and
(\ref{eq:3.5.9}) to the right hand side of (\ref{eq:3.5.2}) gives us
\begin{equation}
\label{eq:3.5.10}
\| d_{A(t)}^* F_{A(t)} \| \leq \| \ddot{a}(t) \| + \| \dot{a}(t)
\| + C \t | \dot{a}(t) |_{C^{k, \mu}} \leq C \| \dot{a} \|_{[t-1, t+1]} 
\end{equation}
We have also from elliptic estimates for Sobolev norms applied to
the equation $ d^*A d_A \be = - *[a, * \dot{a}]$, that
\begin{equation}
\label{eq:3.5.11}
\| d_{A} \be (t) \| \leq C \tau \| \dot{a}(t)\|
\end{equation}
If $\tau$ is small, (\ref{eq:3.5.11}) implies
\begin{equation}
\label{eq:3.5.12}
\frac 12 \| \dot{a}(t) \| \leq \| \eta (t) \| \leq 2 \| \dot{a}(t) \|
\end{equation}
Putting together \eqref{eq:3.5.9}, (\ref{eq:3.5.10}) and
(\ref{eq:3.5.12}) and plugging in both sides (\ref{eq:3.5.1}), we have
\begin{equation}
\label{eq:3.5.13}
\int_T^{\infty} \| \dot{a} \|^2 ds \leq C ( \int_{T-1}^{T+1} \|
\dot{a} \|^2 )^{\frac{1}{2(1-\theta)}} + C \int_{T-1}^{T+1} \|
\dot{a} \|^2  \leq C ( \int_{T-1}^{T+1} \|
\dot{a} \|^2 )^{\frac{1}{2(1-\theta)}}
\end{equation}
where $\theta \in (0, \frac12)$ depend only on $A_0$.
(\ref{eq:3.5.13}) gives an integral decay estimate for
$\|\dot{a}\|$. Recall that $|a(t)|_{C^{k,\mu}} \leq \t$, for $ R < t <
\infty$.  Now it is an easy analytical exercise to show (for example,
as in \cite[3.15]{SL2} ) that there exists $T_1>0$, $\al >0$, such that
\begin{equation}
\label{eq:3.5.14} 
\int_t^{\infty} \|\dot{a}(s) \|ds < C t^{-\al}, \quad \hbox{for } t
\geq T_1
\end{equation}
Therefore
\begin{equation}
\label{eq:3.5.15} 
\| A(t) - A_0 \| \leq C t^{-\al},\quad \hbox{for } t
\geq T_1
\end{equation}
and by elliptic estimates,
\begin{equation}
\label{eq:3.5.16} 
|\tilde{A}(t)- A_0|_{C^{k, \mu}} \leq C(k)t^{-\al}, \quad \hbox{for } t
\geq T_1.
\end{equation}
The desired rate of convergence is obtained and the proof of
Theorem~\ref{th:1.1} is finished.

With the convergence from Theorem~\ref{th:1.1} and the integrablity
assumption, the fast convergence of Theorem~\ref{th:1.2} is a
well-known result (see for example the proof of Theorem 1 (i) in
\cite{AS}). We remark here that a proof of Theorem~\ref{th:1.2}
without using the variational inequalities in Section~\ref{sec:3.3} is
possible. In fact, the variational approach may be totally avoided in
this case as in Cheeger and Tian \cite{CT}, where integrablity of the
cone is assumed.


\section[A result for Yang-Mills flows]{A result of existence and convergence for Yang-Mills flows}

In this section we give an application of the previous methods to
Yang-Mills flows. We shall show that a flow which starts from a
connection sufficiently close (in smooth norms) to a smooth local
minimizer of the Yang-Mills functional will converge asymptotically to
a smooth Yang-Mills connection near the minimizer. Our method, like
before, still consists of two steps, first we choose a suitable gauge,
and then we use the result for parabolic evolution equations (Theorem
2) in \cite{SL1} .

Consider the the following Yang-Mills flow equation for connections on
bundle $E$ on Riemannian manifolds $M$
\begin{equation}
\label{eq:4.1.1} 
\frac{\partial}{\partial t}A(t) = - d_{A(t)}^*F_{A(t)}
\end{equation}
Idealistically, if (\ref{eq:4.1.1}) has a solution $A(t)$ on $[0,
\infty)$, the limit of $A(t)$ at $\infty$ should be a Yang-Mills
connection. Then this will give us a way to homotopically deform an
arbitrary connection into a Yang-Mills connection and hopefully we can
have a Morse theory suitably defined. However, the long-range
existence of solutions of (\ref{eq:4.1.1}) as well as the existence and
regularity of the limit in general are not at all obvious. Nonetheless,
near a local minimizer of the Yang-Mills connection, we are able to
show the flow does exist for all time and converges.

We first note that (\ref{eq:4.1.1}) is not parabolic due to the fact
that $d_A^*F_A$ is not elliptic in $A$. As before, we hope to use the
Coulomb gauge to make the equation parabolic. We note that
(\ref{eq:4.1.1}) actually implies
\begin{equation}
\label{eq:4.1.2} 
d_A^*(\dot{A}) = 0,
\end{equation}
where $\dot{A} = \frac{\partial}{\partial t} A(t)$. This follows from
$d^*_A d^*_A F_A = 0$ by \eqref{eq:3.3.12}.

Assume $A_0$ is a fixed smooth
Yang-Mills connection, we have the following
theorem now. Fix $l$ integer such that $H^{l}(S^{n-1}) \subset C^{3,
\mu}(S^{n-1})$.   
\begin{theorem}
\label{th:4.1} 
There exists $\e = \e(A_0) > 0$, $\al = \al(A_0) >0$ such that for any
given smooth $a_0 \in \Omega^1(\fg_E)$ with $\| a_0 \|_{H^{l+2}} < \e$,
there is a $T_* >0 $ and $A(t)$, a $C^\infty(M \times [0, T^*))$
solution of (\ref{eq:4.1.1}) satisfying $A(0) = A_0 + a_0$, $\sup_{[0, T_*)}
\| A(t) - A_0\|_{H^{l}}  < \e^\al$ and either
\begin{equation}
\label{eq:4.1.1a} 
T_* < \infty \hbox{  and  } \lim_{t \uparrow T_*} \YM(A(t)) \leq
\YM (A_0) - \e
\end{equation}
or
\begin{equation}
\label{eq:4.1.1b} 
T_* = \infty \hbox{  and  } \lim_{t \ra \infty} (|\dot{A}(t)|_{C^1} +
| A(t) - A_1|_{C^2}) = 0
\end{equation}
where $A_1$ is a smooth Yang-Mills connection on $M$.
\end{theorem}

Assume $I$ is an interval of possibly infinite length and $A(t)
+ \be(t) dt$ is a smooth connection on $E \times I$, where $A(t)$ are
connections on $E$ and $\be(t) \in \Gamma(\fg_E)$. We also assume that
under a gauge transformation $g \in \Gamma( \Aut P \times I)$,
\begin{equation}
\label{eq:4.1.3}
g(A(t) + \be(t)dt) = A_1(t)
\end{equation}
where $A_1(t)$ are connections on $E$ and satisfies (\ref{eq:4.1.1})
and hence (\ref{eq:4.1.2}) for $t \in I$. From (\ref{eq:4.1.3}), we
obtain,
\begin{align}
   \label{eq:4.1.4} 
& A_1 = g A g^{-1}  - dg \, g^{-1}\\
\label{eq:4.1.5} 
& \frac{\partial}{\partial t}g = g \be
\end{align}
We observe that (\ref{eq:4.1.4}) implies that $A_1$ and $A$ are gauge
equivalent connections on $M$.  Substitute (\ref{eq:4.1.4}) and
(\ref{eq:4.1.5}) into (\ref{eq:4.1.1}) and (\ref{eq:4.1.2}), by
straightforward computation, we obtain the following equations for $A$
and $\be$
\begin{align}
\label{eq:4.1.6} 
& \dot{A} = -d_A^*F_A + d_A \be\\
\label{eq:4.1.7} 
& d_A^*(\dot{A} - d_A \be) =0
\end{align}
We shall view (\ref{eq:4.1.6}) and (\ref{eq:4.1.7}) (which is implied by
(\ref{eq:4.1.6})) as equations for the connection $A + \be dt$ on $M
\times I$. It is easy to see that they are gauge invariant equations
for gauge transformations on bundle $E \times I$, i.e., if $g \in
\Gamma (\Aut P \times I)$, $g(A +\be
dt) = A_1 + \be_1 dt$ then $A_1 + \be_1 dt$ also satisfies
(\ref{eq:4.1.6}) and (\ref{eq:4.1.7}). This point of view enables us
to consider as before a standard form of $A + \be dt$ around a
connection $A_0$ on $E$ and make the system (\ref{eq:4.1.6}) parabolic
in $A$ and (\ref{eq:4.1.7}) elliptic in $\be$. If $A + \be dt = A_0 +
a(t) + \be(t)dt$ is under such a standard form, i.e. $d_{A_0}^* a =0$
and $\be \in \Ker(d_{A_0})^{\perp}$, then we may solve $\be$ from
(\ref{eq:4.1.7}) and by substituting in (\ref{eq:4.1.6}), rewrite
(\ref{eq:4.1.6}) as
\begin{equation}
\label{eq:4.1.8}
\dot{a} = -d_{A_0 +a}^* F_{A_0 + a} - d_{A_0}d_{A_0}^* a + d_A G_A
(*[a, *\dot{a}])
\end{equation}
where $G_A = (\Delta_A)^{-1}: \Im(d_A^*) \ra \Ker(d_{A_0})^{\perp}$.

From above, we observe that instead of proving Theorem~\ref{th:4.1} in
terms of (\ref{eq:4.1.1}), it suffices to prove the same conclusions
hold for $A(t) = A_0 + a(t)$, $a(t) \in \Ker(d_{A_0}^*)$ with $a(t)$
being solution to (\ref{eq:4.1.8}) with initial value $a_0$ (up to a gauge transformation, we may assume $d_{A_0}^*a_0 = 0$). For if we
prove the latter, a solution to (\ref{eq:4.1.1})  may be obtained by
$$ \bar{A}(t) = g(A(t) + \be(t)),~~ \frac{\partial }{\partial t}g = g\be $$
where $\be = G_A (*[a, *\dot{a}])$. We can show that $g(t) \ra g_0$ for
some $g_0$ and $\dot{g} \ra 0$ in $C^{k}$ as $t \ra \infty$, therefore
$\bar{A}(t)$ will have limit at infinity $g_0(A_1)$ if $A_1$ is the
limit of $A(t)$ and the same conclusions hold for $\bar{A}$.

After this observation, notice that (\ref{eq:4.1.8}) is essentially in
the form of equation (0.1) in \cite{SL1} and the proof there applies
with slight adjustment. We remark that the variation inequalities in
Section 3.3 are again used. Since there is no difficulty, we shall omit
the details here.  We have the following obvious corollary from
Theorem~\ref{th:4.1} .
\begin{corollary}
\label{co:4.2}
If $A_0$ is a smooth local minimizer of Yang-Mills functional on $E$. Then
there exists $\e = \e(A_0) >0$ and $\al = \al(A_0) >0$ such that for any
given smooth $a_0 \in
\Omega^1(\fg_E)$ with with $\| a_0 \|_{H^{l+2}} < \e$,
there is a  $A(t)$, a $C^\infty(M \times [0, \infty))$ solution of
(\ref{eq:4.1.1}) satisfying $A(0) = A_0 + a_0$, and
\begin{equation}
\label{eq:4.1.1.2}
\lim_{t \ra \infty} (|\dot{A}(t)|_{C^1} +
| A(t) - A_1|_{C^2}) = 0
\end{equation}
where $A_1$ is a smooth Yang-Mills connection on $M$ with
$|A_1|_{H^{l}} \leq \e^\al$.
\end{corollary}

{\small

\vspace{2ex}
\noindent {\sc Baozhong Yang\\
Department of Mathematics\\
Stanford University, CA 94305\\}
Email: {\tt byang@math.stanford.edu}
}

\end{document}